\documentclass[thmsa,12pt]{article}
\usepackage{amssymb}
\usepackage{sw20lart}

\input{tcilatex}

\begin{document}

\begin{center}
\bigskip

\textbf{\textit{G}\ Method\ and Finite-Time Consensus}

\bigskip

\textbf{Udrea\ P\u{a}un}

\bigskip
\end{center}

\noindent We give an extension of the $G$ method, with results, the
extension and results being partly suggested by the finite Markov chains and
specially by the finite-time consensus problem for the DeGroot model and
that for the DeGroot model on distributed systems. For the (homogeneous and
nonhomogeneous) DeGroot model, using the $G$ method, a result for reaching a
partial or total consensus in a finite time is given. Further, we consider a
special submodel/case of the DeGroot model, with examples and comments --- a
subset/subgroup property is discovered. For the DeGroot model on distributed
systems, using the $G$ method too, we have a result for reaching a partial
or total (distributed) consensus in a finite time similar to that for the
DeGroot model for reaching a partial or total consensus in a finite time.
Then we show that for any connected graph having $2^{m}$ vertices, $m\geq 1,$
and a spanning subgraph isomorphic to the $m$-cube graph, distributed
averaging is performed in $m$ steps --- this result can be extended ---
research work --- for any graph with $n_{1}n_{2}...n_{t}$ vertices under
certain conditions, where $t,$ $n_{1},n_{2},...,n_{t}\geq 2,$ and, in this
case, distributed averaging is performed in $t$ steps.

\medskip


\medskip


\begin{center}
\bigskip

\textbf{1. AN\ EXTENSION\ OF\ \textit{G}\ METHOD}

\bigskip
\end{center}

In this section, we give an extension of the $G$ method from [18], see
also[17] and [19] --- a good thing both for the forward products of
matrices, in particular, for the finite Markov chains, and for the backwards
products of matrices, in particular, for the DeGroot model and for the
DeGroot model on distributed systems.

\smallskip

Set 
\[
\text{Par}\left( E\right) =\left\{ \Delta \left| \text{ }\Delta \text{ is a
partition of }E\right. \right\} , 
\]

\noindent where $E$ is a nonempty set. We shall agree that the partitions do
not contain the empty set. $\left( E\right) $ is the improper (degenerate)
partition of $E$ --- we need this partition.

\smallskip

\textbf{Definition\ 1.1.} Let $\Delta _{1},\Delta _{2}\in $Par$\left(
E\right) .$ We say that $\Delta _{1}$\textit{\ is finer than }$\Delta _{2}$
if $\forall V\in \Delta _{1},$ $\exists W\in \Delta _{2}$ such that $%
V\subseteq W.$

Write $\Delta _{1}\preceq \Delta _{2}$ when $\Delta _{1}$ is finer than $%
\Delta _{2}.$

\smallskip

In this article, a vector is a row vector and a stochastic matrix is a row
stochastic matrix.

\smallskip

Set 
\[
\left\langle m\right\rangle =\left\{ 1,2,...,m\right\} \text{ (}m\geq 1\text{%
)}. 
\]

\smallskip

Let $Z$ be an $m\times n$ matrix. $Z_{ij}$ or, if confusion can arise, $%
Z_{i\rightarrow j}$ denotes the entry $\left( i,j\right) $ of matrix $Z,$ $%
\forall i\in \left\langle m\right\rangle ,$ $\forall j\in \left\langle
n\right\rangle .$

\smallskip

Set 
\[
C_{m,n}=\left\{ P\left| \text{ }P\text{ is a complex }m\times n\text{ matrix}%
\right. \right\} \text{,} 
\]
\[
R_{m,n}=\left\{ P\left| \text{ }P\text{ is a real }m\times n\text{ matrix}%
\right. \right\} \text{,} 
\]
\[
N_{m,n}=\left\{ P\left| \text{ }P\text{ is a nonnegative }m\times n\text{
matrix}\right. \right\} \text{,} 
\]
\[
S_{m,n}=\left\{ P\left| \text{ }P\text{ is a stochastic }m\times n\text{
matrix}\right. \right\} 
\]
\noindent (... the stochastic matrices are defined in Definition 1.3), 
\[
C_{n}=C_{n,n}\text{,} 
\]
\[
R_{n}=R_{n,n}\text{,} 
\]
\[
N_{n}=N_{n,n}\text{,} 
\]
\[
S_{n}=S_{n,n}. 
\]

Let $P=\left( P_{ij}\right) \in C_{m,n}$. Let $\emptyset \neq U\subseteq
\left\langle m\right\rangle $ and $\emptyset \neq V\subseteq \left\langle
n\right\rangle $. Set the matrices 
\[
P_{U}=\left( P_{ij}\right) _{i\in U,j\in \left\langle n\right\rangle },\text{
}P^{V}=\left( P_{ij}\right) _{i\in \left\langle m\right\rangle ,j\in V},%
\text{ and }P_{U}^{V}=\left( P_{ij}\right) _{i\in U,j\in V}. 
\]

\smallskip

\textbf{Remark 1.2.} (Basic properties of the operators $\left( \cdot
\right) _{\left( \cdot \right) },\left( \cdot \right) ^{\left( \cdot \right)
},$ and $\left( \cdot \right) _{\left( \cdot \right) }^{\left( \cdot \right)
}$ on the products of complex matrices.) Let $P_{1}\in C_{n_{1},n_{2}},$ $%
P_{2}\in C_{n_{2},n_{3}},$ $...,$ $P_{t}\in C_{n_{t},n_{t+1}},$ where $t\geq
2$. Let $k\in \left\langle t-1\right\rangle .$ Let $\emptyset \neq
U\subseteq \left\langle n_{1}\right\rangle $ and $\emptyset \neq V\subseteq
\left\langle n_{t+1}\right\rangle $. Obviously, we have

(a) 
\[
\left( P_{1}P_{2}...P_{t}\right) _{U}=\left( P_{1}P_{2}...P_{k}\right)
_{U}P_{k+1}P_{k+2}...P_{t}; 
\]

(b) 
\[
\left( P_{1}P_{2}...P_{t}\right) ^{V}=P_{1}P_{2}...P_{k}\left(
P_{k+1}P_{k+2}...P_{t}\right) ^{V}; 
\]

(c) 
\[
\left( P_{1}P_{2}...P_{t}\right) _{U}^{V}=\left( P_{1}P_{2}...P_{k}\right)
_{U}\left( P_{k+1}P_{k+2}...P_{t}\right) ^{V}; 
\]

(d) 
\[
\left( P_{1}P_{2}...P_{t}\right) _{U}^{V}=\sum\limits_{W\in \Delta }\left(
P_{1}P_{2}...P_{k}\right) _{U}^{W}\left( P_{k+1}P_{k+2}...P_{t}\right)
_{W}^{V},\text{ }\forall \Delta \in \text{Par}\left( \left\langle
n_{k+1}\right\rangle \right) ; 
\]

(e) 
\[
\left( P_{1}P_{2}...P_{t}\right) _{U}^{V}\geq \left(
P_{1}P_{2}...P_{k}\right) _{U}^{W}\left( P_{k+1}P_{k+2}...P_{t}\right)
_{W}^{V},\text{ }\forall W,\text{ }\emptyset \neq W\subseteq \left\langle
n_{k+1}\right\rangle . 
\]
\noindent ((c)$\Longrightarrow $(d)$\Longrightarrow $(e))

\smallskip

Let $W$ be a nonempty finite set. Suppose that $W=\left\{
s_{1},s_{2},...,s_{t}\right\} $. Set 
\[
\left( \left\{ i\right\} \right) _{i\in W}\in \text{Par}\left( W\right) ,%
\text{ }\left( \left\{ i\right\} \right) _{i\in W}=\left( \left\{
s_{1}\right\} ,\left\{ s_{2}\right\} ,...,\left\{ s_{t}\right\} \right) . 
\]

\noindent \textit{E.g.}, 
\[
\left( \left\{ i\right\} \right) _{i\in \left\langle n\right\rangle }=\left(
\left\{ 1\right\} ,\left\{ 2\right\} ,...,\left\{ n\right\} \right) . 
\]

\noindent $\left( \left\{ i\right\} \right) _{i\in W}$ is the finest
partition of $W$ while $\left( W\right) $ is the coarsest partition of $W$.
Let $L\subset W$. Let $K=L^{c}$, $L^{c}$ is the complement of $L$ --- so, $%
\emptyset \neq K\subseteq W.$ Suppose that $K=\left\{
s_{j_{1}},s_{j_{2}},...,s_{j_{u}}\right\} $ ($1\leq u\leq t$). Set 
\[
\left( \left\{ i\right\} \!,L\right) _{i\in K}\in \text{Par}\left( W\right) ,%
\text{ }\left( \left\{ i\right\} \!,L\right) _{i\in K}=\left\{ 
\begin{array}{l}
\left( \left\{ i\right\} \right) _{i\in W}\text{ if }\left| L\right|
=0,\smallskip \\ 
\left( \left\{ s_{j_{1}}\right\} \!,\left\{ s_{j_{2}}\right\}
\!,...,\!\left\{ s_{j_{u}}\right\} \!,L\right) \text{ if }\left| L\right|
\geq 1%
\end{array}
\right. 
\]
\noindent ($\left| \cdot \right| $ is the cardinality). Obviously, $\left(
\left\{ i\right\} \!,L\right) _{i\in K}=\left( \left\{ i\right\} \right)
_{i\in W}$ if $\left| L\right| =1$. So, $\left( \left\{ i\right\}
\!,L\right) _{i\in K}=\left( \left\{ i\right\} \right) _{i\in W}$ if $\left|
L\right| \leq 1$. Further, we suppose that $t\geq 2$. Let $\left(
L_{1},L_{2},...,L_{m}\right) \in $Par$\left( W\right) $, where $m\geq 2$ ($%
L_{1},L_{2},...,L_{m}\neq \emptyset $; $m\geq 2\Longrightarrow t\geq 2$).
Suppose that $L_{1}=\left\{ s_{l_{1}},s_{l_{2}},...,s_{l_{v}}\right\} $ ($%
v\geq 1$). Set 
\[
\left( \left\{ i\right\} ,L_{2},...,L_{m}\right) _{i\in L_{1}}\in \text{Par}%
\left( W\right) , 
\]
\[
\left( \left\{ i\right\} ,L_{2},...,L_{m}\right) _{i\in L_{1}}=\left(
\left\{ s_{l_{1}}\right\} ,\left\{ s_{l_{2}}\right\} ,...,\left\{
s_{l_{v}}\right\} ,L_{2},...,L_{m}\right) . 
\]

\noindent \textit{E.g.}, 
\[
\left( \left\{ i\right\} ,\left\{ 4,5\right\} ,\left\{ 6,7,8\right\} \right)
_{i\in \left\langle 3\right\rangle }=\left( \left\{ 1\right\} ,\left\{
2\right\} ,\left\{ 3\right\} ,\left\{ 4,5\right\} ,\left\{ 6,7,8\right\}
\right) . 
\]
\noindent $\left( \left\{ i\right\} ,L\right) _{i\in K}$ is the finest
partition of $W$ when $\left| L\right| \leq 1$; $\left( \left\{ i\right\}
,L\right) _{i\in K}$ is the finest partition (of $W$) among the partitions
of $W$ which contain $L$ when $\left| L\right| >1$ while $\left( \left\{
i\right\} ,L_{2},...,L_{m}\right) _{i\in L_{1}}$ is the finest partition (of 
$W$) among the partitions of $W$ which contain $L_{2},...,L_{m}$. Obviously, 
$\left( \left\{ i\right\} ,L_{2},...,L_{m}\right) _{i\in L_{1}}=\left(
\left\{ i\right\} \right) _{i\in W}$ when $\left| L_{2}\right| =...=\left|
L_{m}\right| =1.$ Obviously, the partition $\left( \left\{ i\right\}
,L\right) _{i\in K}$ when $\left| L\right| \geq 1$ is a special one of $%
\left( \left\{ i\right\} ,L_{2},...,L_{m}\right) _{i\in L_{1}}.$

\smallskip

\textbf{Definition 1.3.} Let $P\in C_{m,n}$. We say that $P$ is a \textit{%
generalized stochastic }(\textit{complex})\textit{\ matrix} if $\exists a\in 
\mathbb{C}$ ($a=0$ or $a\neq 0$) such that 
\[
\sum_{j\in \left\langle n\right\rangle }P_{ij}=a,\text{ }\forall i\in
\left\langle m\right\rangle . 
\]

\noindent (We have generalized stochastic complex matrices, generalized
stochastic real matrices, generalized stochastic nonnegative matrices... If $%
P\in N_{m,n},$ then $P$ is a generalized stochastic (nonnegative) matrix if
and only if $\exists a\geq 0,$ $\exists Q\in S_{m,n}$ such that $P=aQ$ ---
see, \textit{e.g}., [17, Definition 1.2].) The generalized stochastic
matrices with row sums equal to 1 are called \textit{stochastic complex }%
when the matrices are complex, \textit{stochastic real }when the matrices
are real, \textit{stochastic nonnegative} or, for short, \textit{stochastic }%
when the matrices are nonnegative...

\smallskip

\textbf{Definition 1.4.} (A generalization of Definition 1.3 from [17].) Let 
$P\in C_{m,n}$. Let $\Delta \in $Par$\left( \left\langle m\right\rangle
\right) $ and $\Sigma \in $Par$\left( \left\langle n\right\rangle \right) $.
We say that $P$ is a $\left[ \Delta \right] $\textit{-stable matrix on }$%
\Sigma $ if $P_{K}^{L}$ is a generalized stochastic matrix, $\forall K\in
\Delta ,$ $\forall L\in \Sigma $. In particular, a $\left[ \Delta \right] $%
-stable matrix on $\left( \left\{ j\right\} \right) _{j\in \left\langle
n\right\rangle }$ is called $\left[ \Delta \right] $\textit{-stable }for
short.

\smallskip

\textbf{Definition 1.5.} (A generalization of Definition 1.4 from [17].) Let 
$P\in C_{m,n}$. Let $\Delta \in $Par$\left( \left\langle m\right\rangle
\right) $ and $\Sigma \in $Par$\left( \left\langle n\right\rangle \right) $.
We say that $P$ is a $\Delta $\textit{-stable matrix on }$\Sigma $ if $%
\Delta $ is the least fine partition for which $P$ is a $\left[ \Delta %
\right] $-stable matrix on $\Sigma $. In particular, a $\Delta $-stable
matrix on $\left( \left\{ j\right\} \right) _{j\in \left\langle
n\right\rangle }$ is called $\Delta $\textit{-stable} while a $\left(
\left\langle m\right\rangle \right) $-stable matrix on $\Sigma $ is called 
\textit{stable on }$\Sigma $ for short. A stable matrix on $\left( \left\{
j\right\} \right) _{j\in \left\langle n\right\rangle }$ is called \textit{%
stable} for short.

\smallskip

For the $\left[ \Delta \right] $- and $\Delta $- stable matrices on $\Sigma $%
, $\Delta \in $Par$\left( \left\langle m\right\rangle \right) $, $\Sigma \in 
$Par$\left( \left\langle n\right\rangle \right) $, we will use the generic
name ``generalized stable (complex) matrices'' --- we have generalized
stable complex matrices, generalized stable real matrices, generalized
stable nonnegative matrices, generalized stable stochastic matrices...

\smallskip

The stable matrices have a central role in the finite Markov chain theory,
in the DeGroot model theory (see Section 2), in the theory of DeGroot model
on distributed systems (see Section 3), ... The stable matrices on $\Sigma $
are important for the finite Markov chains (see, \textit{e.g.}, [17, p.
389]), for the partial consensuses (see, in Section 2, Theorem 2.5; see also
Example 2.8), for the partial distributed consensuses (see, in Section 3,
Remark 3.2(b)), ...

\smallskip

\textbf{Example 1.6.} Let 
\[
P=\left( 
\begin{array}{ccccc}
2\medskip & 3 & \frac{1}{2} & \frac{1}{2} & 0 \\ 
1 & 4\medskip & 0 & \frac{1}{2} & \frac{1}{2} \\ 
4 & 6 & 0\medskip & 0 & 7 \\ 
9 & 1 & 1 & 2 & 4%
\end{array}
\right) ,\text{ }Q=\left( 
\begin{array}{ccc}
\frac{1}{2}\medskip & 0 & 0 \\ 
\frac{1}{2} & 0\medskip & 0 \\ 
1 & 2 & 3\medskip \\ 
1 & 2 & 3%
\end{array}
\right) , 
\]
\[
T=\left( 
\begin{array}{cccc}
3\medskip & -\frac{1}{4} & \frac{1}{4} & \frac{2}{4} \\ 
3 & -\frac{1}{4} & \frac{3}{4} & 0%
\end{array}
\right) ,\text{ and }Z=\left( 
\begin{array}{cc}
\frac{1}{4}\medskip & -\frac{3}{4} \\ 
\frac{1}{4} & -\frac{3}{4}%
\end{array}
\right) . 
\]

\noindent $P$ is $\left[ \left( \left\{ 1,2\right\} ,\left\{ 3,4\right\}
\right) \right] $-stable on $\left( \left\{ 1,2\right\} ,\left\{
3,4,5\right\} \right) $ because $P_{K}^{L}$ is a generalized stochastic
matrix, $\forall K\in \left( \left\{ 1,2\right\} ,\left\{ 3,4\right\}
\right) ,$ $\forall L\in \left( \left\{ 1,2\right\} ,\left\{ 3,4,5\right\}
\right) $ --- indeed, 
\[
P_{\left\{ 1,2\right\} }^{\left\{ 1,2\right\} }=\left( 
\begin{array}{cc}
2 & 3 \\ 
1 & 4%
\end{array}
\right) , 
\]
\noindent and it is a generalized stochastic matrix because $2+3=1+4,$%
\[
P_{\left\{ 1,2\right\} }^{\left\{ 3,4,5\right\} }=\left( 
\begin{array}{ccc}
\frac{1}{2}\medskip & \frac{1}{2} & 0 \\ 
0 & \frac{1}{2} & \frac{1}{2}%
\end{array}
\right) , 
\]

\noindent and it is a generalized stochastic matrix, etc. Moreover, $P$ is $%
\left[ \Delta \right] $-stable on $\left( \left\{ 1,2\right\} ,\left\{
3,4,5\right\} \right) ,$ $\forall \Delta \preceq \left( \left\{ 1,2\right\}
,\left\{ 3,4\right\} \right) .$ Moreover, $P$ is $\left( \left\{ 1,2\right\}
,\left\{ 3,4\right\} \right) $-stable on $\left( \left\{ 1,2\right\}
,\left\{ 3,4,5\right\} \right) .$ $Q$ is $\left[ \left( \left\{ 1,2\right\}
,\left\{ 3,4\right\} \right) \right] $-stable 
\[
\text{(}Q_{\left\{ 1,2\right\} }^{\left\{ 1\right\} }=\left( 
\begin{array}{c}
\frac{1}{2}\medskip \\ 
\frac{1}{2}%
\end{array}
\right) =\frac{1}{2}\left( 
\begin{array}{c}
1\medskip \\ 
1%
\end{array}
\right) ,\text{ }Q_{\left\{ 1,2\right\} }^{\left\{ 2\right\} }\equiv
Q_{\left\{ 1,2\right\} }^{\left\{ 3\right\} }=\left( 
\begin{array}{c}
0 \\ 
0%
\end{array}
\right) ,\text{ ...} 
\]

\noindent --- $\equiv $ is the identical sign). Moreover, it is $\left[
\Delta \right] $-stable, $\forall \Delta \preceq \left( \left\{ 1,2\right\}
,\left\{ 3,4\right\} \right) .$ Moreover, it is $\left( \left\{ 1,2\right\}
,\left\{ 3,4\right\} \right) $-stable. $T$ is stable on $\left( \left\{
1\right\} ,\left\{ 2\right\} ,\left\{ 3,4\right\} \right) $ --- it follows
that $T^{\left\langle 2\right\rangle }$ is stable (the rows of $%
T^{\left\langle 2\right\rangle }$ are identical --- any stable matrix has
the rows identically). $Z$ is stable (the rows of $Z$ are identical).

\smallskip

Let $\Delta \in $Par$\left( \left\langle m\right\rangle \right) $ and $%
\Sigma \in $Par$\left( \left\langle n\right\rangle \right) $. Set (see [17]
for $G_{\Delta ,\Sigma }$ and [18] for $\overline{G}_{\Delta ,\Sigma }$) 
\[
G_{\Delta ,\Sigma }=G_{\Delta ,\Sigma }\left( m,n\right) =\left\{ P\left| 
\text{ }P\in S_{m,n}\text{ and }P\text{ is a }\left[ \Delta \right] \text{%
-stable matrix on }\Sigma \right. \right\} , 
\]
\[
\overline{G}_{\Delta ,\Sigma }=\overline{G}_{\Delta ,\Sigma }\left(
m,n\right) =\left\{ P\left| \text{ }P\in N_{m,n}\text{ and }P\text{ is a }%
\left[ \Delta \right] \text{-stable matrix on }\Sigma \right. \right\} , 
\]
\[
\widetilde{G}_{\Delta ,\Sigma }=\widetilde{G}_{\Delta ,\Sigma }\left(
m,n\right) =\left\{ P\left| \text{ }P\in R_{m,n}\text{ and }P\text{ is a }%
\left[ \Delta \right] \text{-stable matrix on }\Sigma \right. \right\} , 
\]

\noindent and 
\[
\widehat{G}_{\Delta ,\Sigma }=\widehat{G}_{\Delta ,\Sigma }\left( m,n\right)
=\left\{ P\left| \text{ }P\in C_{m,n}\text{ and }P\text{ is a }\left[ \Delta %
\right] \text{-stable matrix on }\Sigma \right. \right\} 
\]
\noindent ($G_{\Delta ,\Sigma }$ is a set of generalized stable stochastic
matrices, $\overline{G}_{\Delta ,\Sigma }$ is a set of generalized stable
nonnegative matrices, ...).

\smallskip

When we study or even when we construct products of complex matrices (in
particular, products of stochastic matrices) using $G_{\Delta ,\Sigma },$ $%
\overline{G}_{\Delta ,\Sigma },$ $\widetilde{G}_{\Delta ,\Sigma },$ or $%
\widehat{G}_{\Delta ,\Sigma },$ we shall refer this as the $G$\textit{\
method}. $G$ comes from the verb \textit{to group} and its derivatives.

\smallskip

The $G$\ method can more be extended considering other sets of matrices
besides $S_{m,n},$ $N_{m,n},$ $R_{m,n},$ and $C_{m,n},$ then, if necessary,
generalizing Definition 1.3...

\smallskip

The $G$\ method has serious applications, see, \textit{e.g.}, [20]-[21].

\smallskip

Since $G_{\Delta ,\Sigma }\subseteq \overline{G}_{\Delta ,\Sigma }\subseteq 
\widetilde{G}_{\Delta ,\Sigma }\subseteq \widehat{G}_{\Delta ,\Sigma }$ ($%
\Delta \in $Par$\left( \left\langle m\right\rangle \right) ,$ $\Sigma \in $%
Par$\left( \left\langle n\right\rangle \right) $), any property which holds
for $\widehat{G}_{\Delta ,\Sigma }$ also holds for $G_{\Delta ,\Sigma },$ $%
\overline{G}_{\Delta ,\Sigma },$ and $\widetilde{G}_{\Delta ,\Sigma },$ any
property which holds for $\widetilde{G}_{\Delta ,\Sigma }$ also holds for $%
G_{\Delta ,\Sigma }$ and $\overline{G}_{\Delta ,\Sigma }$... Several
properties of $G_{\Delta ,\Sigma }$ were given in Remark 2.1 from [17] ---
some of these properties can be extended for $\overline{G}_{\Delta ,\Sigma
}, $ for $\widetilde{G}_{\Delta ,\Sigma },$ or for $\widehat{G}_{\Delta
,\Sigma }.$ Below we give a few properties of $G_{\Delta ,\Sigma },$ $%
\overline{G}_{\Delta ,\Sigma },$ $\widetilde{G}_{\Delta ,\Sigma },$ or $%
\widehat{G}_{\Delta ,\Sigma }.$

\smallskip

\textbf{Remark 1.7.} (a) 
\[
S_{m,n}=G_{\left( \left\langle m\right\rangle \right) ,\left( \left\langle
n\right\rangle \right) }=G_{\left( \left\{ i\right\} \right) _{i\in
\left\langle m\right\rangle },\Sigma },\text{ }\forall m,n\geq 1,\text{ }%
\forall \Sigma \in \text{Par}\left( \left\langle n\right\rangle \right) 
\]
\noindent ($\left( \left\langle m\right\rangle \right) $ is the improper
(degenerate) partition of $\left\langle m\right\rangle $...; $\Sigma $ can
be equal to $\left( \left\{ j\right\} \right) _{j\in \left\langle
n\right\rangle }$ (so, $S_{m,n}=G_{\left( \left\{ i\right\} \right) _{i\in
\left\langle m\right\rangle },\left( \left\{ j\right\} \right) _{j\in
\left\langle n\right\rangle }}$), or, more generally, to $\left( \left\{
j\right\} ,K^{c}\right) _{j\in K},$ $\emptyset \neq K\subseteq \left\langle
n\right\rangle ,$ or to...).

\smallskip

(b) 
\[
N_{m,n}\supset \overline{G}_{\left( \left\langle m\right\rangle \right)
,\left( \left\langle n\right\rangle \right) },\text{ }R_{m,n}\supset 
\widetilde{G}_{\left( \left\langle m\right\rangle \right) ,\left(
\left\langle n\right\rangle \right) },\text{ }C_{m,n}\supset \widehat{G}%
_{\left( \left\langle m\right\rangle \right) ,\left( \left\langle
n\right\rangle \right) },\text{ }\forall m\geq 2,\text{ }\forall n\geq 1, 
\]

\noindent and $\overline{G}_{\left( \left\langle m\right\rangle \right)
,\left( \left\langle n\right\rangle \right) }$ is equal to the set of
generalized stochastic nonnegative $m\times n$ matrices, $\widetilde{G}%
_{\left( \left\langle m\right\rangle \right) ,\left( \left\langle
n\right\rangle \right) }$ is equal to the set of generalized stochastic real 
$m\times n$ matrices, and $\widehat{G}_{\left( \left\langle m\right\rangle
\right) ,\left( \left\langle n\right\rangle \right) }$ is equal to the set
of generalized stochastic complex $m\times n$ matrices; 
\[
N_{m,n}=\overline{G}_{\left( \left\{ i\right\} \right) _{i\in \left\langle
m\right\rangle },\Sigma },\text{ }R_{m,n}=\widetilde{G}_{\left( \left\{
i\right\} \right) _{i\in \left\langle m\right\rangle },\Sigma },\text{ }%
C_{m,n}=\widehat{G}_{\left( \left\{ i\right\} \right) _{i\in \left\langle
m\right\rangle },\Sigma }, 
\]
\noindent $\forall m\geq 1,$ $\forall n\geq 1,$ $\forall \Sigma \in $Par$%
\left( \left\langle n\right\rangle \right) .$

\smallskip

(c) The stable complex $m\times n$ matrices on $\Sigma $ are generalized
stochastic, $\forall \Sigma \in $Par$\left( \left\langle n\right\rangle
\right) $ ($\overline{G}_{\left( \left\langle m\right\rangle \right) ,\Sigma
},\widetilde{G}_{\left( \left\langle m\right\rangle \right) ,\Sigma },$ and $%
\widehat{G}_{\left( \left\langle m\right\rangle \right) ,\Sigma }$ are
included in the set of generalized stochastic nonnegative $m\times n$
matrices, that of generalized stochastic real $m\times n$ matrices, and that
of generalized stochastic complex $m\times n$ matrices, respectively).

\smallskip

(d) 
\[
\overline{G}_{\left( \left\langle m\right\rangle \right) ,\left( \left\{
j\right\} ,L_{2},L_{3},...,L_{v}\right) _{j\in L_{1}}}\subseteq \overline{G}%
_{\left( \left\langle m\right\rangle \right) ,\left( \left\{ j\right\}
,\bigcup\limits_{b=2}^{v}L_{b}\right) _{j\in L_{1}}}, 
\]

\noindent $\forall m,$ $n\geq 1,$ $\forall \!\left(
L_{1},\!L_{2},\!...,\!L_{v}\right) \!\in $Par$\left( \left\langle
n\right\rangle \right) \!,$ $v\geq 2$ (obviously, $\left( \left\{ j\right\}
\!\!,L_{2},L_{3},...,L_{v}\right) _{j\in L_{1}}$

\noindent $\in $Par$\left( \left\langle n\right\rangle \right) $). Similarly
for $\widetilde{G}_{\left( \left\langle m\right\rangle \right) ,\left(
\left\{ j\right\} ,L_{2},L_{3},...,L_{v}\right) _{j\in L_{1}}}$and for $%
\widehat{G}_{\left( \left\langle m\right\rangle \right) ,\left( \left\{
j\right\} ,L_{2},L_{3},...,L_{v}\right) _{j\in L_{1}}}$.

\smallskip

Let $P\in \widehat{G}_{\Delta ,\Sigma }.$ Then (see Definition 1.4) $\forall
K\in \Delta ,$ $\forall L\in \Sigma ,$ $\exists a_{K,L}=a_{K,L}\left(
P\right) \in \mathbb{C}$ such that 
\[
\sum\limits_{j\in L}P_{ij}=a_{K,L},\forall i\in K. 
\]
\noindent Set the matrix (for the case when $P\in \overline{G}_{\Delta
,\Sigma },$ see, \textit{e.g.}, also [18]) 
\[
P^{-+}=\left( P_{KL}^{-+}\right) _{K\in \Delta ,L\in \Sigma },\text{ }%
P_{KL}^{-+}=a_{K,L},\text{ }\forall K\in \Delta ,\text{ }\forall L\in \Sigma
. 
\]

\noindent We call the matrix $P^{-+}$ the $\left( \Delta ,\Sigma \right) $-%
\textit{grouped matrix} (\textit{of} $P$), or, for short, if no confusion
can arise, the \textit{grouped matrix} (\textit{of} $P$) --- in [17] this
matrix was called differently... The labels of rows and columns of $P^{-+}$
are $K,$ $K\in \Delta ,$ and $L,$ $L\in \Sigma $, respectively; $%
P_{KL}^{-+}, $ $K\in \Delta ,$ $L\in \Sigma ,$ are the entries of matrix $%
P^{-+}$; $P_{KL}^{-+}$ is the entry $\left( K,L\right) $ of matrix $P^{-+},$ 
$\forall K\in \Delta ,$ $\forall L\in \Sigma .$ If confusion can arise, we
write $P^{-+\left( \Delta ,\Sigma \right) }$ instead of $P^{-+}.$

\smallskip

\textbf{Example 1.8.} Consider the matrices $P,$ $Q,$ $T,$ and $Z$ from
Example 1.6. We have 
\[
P^{-+}=\left( 
\begin{array}{cc}
5 & 1 \\ 
10 & 7%
\end{array}
\right) \text{ (}\Delta =\left( \left\{ 1,2\right\} ,\left\{ 3,4\right\}
\right) ,\text{ }\Sigma =\left( \left\{ 1,2\right\} ,\left\{ 3,4,5\right\}
\right) \text{),} 
\]
\[
Q^{-+}=\left( 
\begin{array}{ccc}
\frac{1}{2} & 0 & 0 \\ 
1 & 2 & 3%
\end{array}
\right) \text{ (}\Delta =\left( \left\{ 1,2\right\} ,\left\{ 3,4\right\}
\right) ,\text{ }\Sigma =\left( \left\{ j\right\} \right) _{j\in
\left\langle 3\right\rangle }\text{),} 
\]
\[
T^{-+}=\left( 
\begin{array}{ccc}
3 & -\frac{1}{4} & \frac{3}{4}%
\end{array}
\right) ,\text{ }Z^{-+}=\left( 
\begin{array}{cc}
\frac{1}{4} & -\frac{3}{4}%
\end{array}
\right) . 
\]

\noindent $T^{-+}$ and $Z^{-+}$ have just one row because $T$ and $Z$ are
stable on $\left( \left\{ 1\right\} ,\left\{ 2\right\} ,\left\{ 3,4\right\}
\right) $ and stable, respectively. Moreover, the row of $\left(
T^{\left\langle 2\right\rangle }\right) ^{-+}=\left( T^{\left\langle
2\right\rangle }\right) ^{-+\left( \left( \left\langle 2\right\rangle
\right) ,\left( \left\{ 1\right\} ,\left\{ 2\right\} \right) \right) }$ is
identical with any row of $T^{\left\langle 2\right\rangle }$ and the row of $%
Z^{-+}$ is identical with any row of $Z.$ Generally speaking, letting $A$
and $B,$ $A,$ $B\in C_{m,n},$ be a stable matrix on $\left( \left\{
j\right\} ,K^{c}\right) _{j\in K},$ $\emptyset \neq K\subseteq \left\langle
n\right\rangle ,$ $\left| K^{c}\right| >1,$ and a stable matrix,
respectively, $\left( A^{K}\right) ^{-+}=\left( A^{K}\right) ^{-+\left(
\left( \left\langle m\right\rangle \right) ,\left( \left\{ j\right\} \right)
_{j\in K}\right) }$ has just one row and this is identical with any row of $%
A^{K}$ and $B^{-+}=B^{-+\left( \left( \left\langle m\right\rangle \right)
,\left( \left\{ j\right\} \right) _{j\in \left\langle n\right\rangle
}\right) }$ has just one row and this is identical with any row of $B$ ---
do not forget! (see Theorem 1.10...)

\smallskip

Below we give a fundamental result.

\smallskip

\textbf{THEOREM\ 1.9.} (i) (a generalization of Theorem 1.5 from [18]) 
\textit{Let }$P_{1}\in \widehat{G}_{\Delta _{1},\Delta _{2}}\subseteq
C_{n_{1},n_{2}}$ \textit{and} $P_{2}\in \widehat{G}_{\Delta _{2},\Delta
_{3}}\subseteq C_{n_{2},n_{3}}.$\textit{\ Then} 
\[
P_{1}P_{2}\in \widehat{G}_{\Delta _{1},\Delta _{3}}\subseteq C_{n_{1},n_{3}}%
\hspace{0.3cm}\text{\textit{and} }\left( P_{1}P_{2}\right) ^{-+}=\text{ }%
P_{1}^{-+}\text{ }P_{2}^{-+}. 
\]

(ii) (a generalization of (i)) \textit{Let }$P_{1}\in \widehat{G}_{\Delta
_{1},\Delta _{2}}\subseteq C_{n_{1},n_{2}},$ $P_{2}\in \widehat{G}_{\Delta
_{2},\Delta _{3}}\subseteq C_{n_{2},n_{3}},$ $...,$ $P_{t}\in \widehat{G}%
_{\Delta _{t},\Delta _{t+1}}\subseteq C_{n_{t},n_{t+1}}$, \textit{where} $%
t\geq 2.$ \textit{Then} 
\[
P_{1}P_{2}...P_{t}\in \widehat{G}_{\Delta _{1},\Delta _{t+1}}\subseteq
C_{n_{1},n_{t+1}}\hspace{0.3cm}\text{\textit{and} }\left(
P_{1}P_{2}...P_{t}\right) ^{-+}=\text{ }P_{1}^{-+}\text{ }%
P_{2}^{-+}...P_{t}^{-+}. 
\]

\smallskip

\textbf{Proof.} (i) Similar to the proof of Theorem 1.5 from [18]. Another
proof, a simple one too, of the first part of (i) is as follows. Let $U\in
\Delta _{1}$ and $V\in \Delta _{3}.$ By Remark 1.2(d) we have 
\[
\left( P_{1}P_{2}\right) _{U}^{V}=\sum_{W\in \Delta _{2}}\left( P_{1}\right)
_{U}^{W}\left( P_{2}\right) _{W}^{V}. 
\]
\noindent Since $\left( P_{1}\right) _{U}^{W}$ and $\left( P_{2}\right)
_{W}^{V}$ are generalized stochastic matrices, $\forall W\in \Delta _{2},$
it follows that 
\[
\left( P_{1}\right) _{U}^{W}\left( P_{2}\right) _{W}^{V}\text{ and }%
\sum_{W\in \Delta _{2}}\left( P_{1}\right) _{U}^{W}\left( P_{2}\right)
_{W}^{V} 
\]
\noindent are generalized stochastic matrices. So, $\left( P_{1}P_{2}\right)
_{U}^{V}$ is a generalized stochastic matrix. Therefore, $P_{1}P_{2}\in 
\widehat{G}_{\Delta _{1},\Delta _{3}}.$

(ii) Induction. 
\endproof%

\smallskip

Recall that in this article, a vector is a row vector. Set $e=e\left(
n\right) =\left( 1,1,...,1\right) \in \mathbb{R}^{n}$ ($e\in \mathbb{C}^{n}$
too), $\forall n\geq 1$; $e^{\prime }$ is its transpose.

\smallskip

Below we give the best result of this section.

\smallskip

\textbf{THEOREM\ 1.10.} (An extension of Theorem 1.6 from [18].) \textit{Let 
}$P_{1}\in \widehat{G}_{\Delta _{1},\Delta _{2}}\subseteq C_{n_{1},n_{2}},$ $%
P_{2}\in \widehat{G}_{\Delta _{2},\Delta _{3}}\subseteq C_{n_{2},n_{3}},$ $%
...,$ $P_{t}\in \widehat{G}_{\Delta _{t},\Delta _{t+1}}\subseteq
C_{n_{t},n_{t+1}}$, \textit{where} $t\geq 2.$\textit{\ Suppose that} 
\[
\Delta _{1}=\left( \left\langle n_{1}\right\rangle \right) \hspace{0.3cm}%
\text{\textit{and }}\Delta _{t+1}=\left( \left\{ j\right\} ,K^{c}\right)
_{j\in K} 
\]
($\emptyset \neq K\subseteq \left\langle n_{t+1}\right\rangle $;\textit{\ }$%
\left| K^{c}\right| \leq 1$ \textit{or} $\left| K^{c}\right| >1$). \textit{%
Then} 
\[
\left( P_{1}P_{2}...P_{t}\right) ^{K} 
\]
\[
\text{(\textit{equivalently}, }P_{1}P_{2}...P_{t-1}P_{t}^{K}\text{ (\textit{%
see Remark} 1.2(b)))} 
\]
\noindent \textit{is a stable matrix} (\textit{i.e., a matrix with identical
rows, see Definition} 1.5 \textit{and Example} 1.6; $\left(
P_{1}P_{2}...P_{t}\right) ^{K}$ \textit{is a submatrix of }$%
P_{1}P_{2}...P_{t},$\textit{\ see the definition of operator }$\left( \cdot
\right) ^{\left( \cdot \right) }$\textit{\ again}; \textit{if }$%
K=\left\langle n_{t+1}\right\rangle ,$\textit{\ then }$\left(
P_{1}P_{2}...P_{t}\right) ^{K}=P_{1}P_{2}...P_{t},$ and $P_{1}P_{2}...P_{t}$%
\textit{\ is a stable matrix}), 
\[
\left( P_{1}P_{2}...P_{t}\right) _{\left\{ i\right\} }^{K}=\left(
P_{1}^{-+}P_{2}^{-+}...P_{t}^{-+}\right) ^{\bigcup\limits_{j\in K}\left\{
\left\{ j\right\} \right\} }= 
\]
\[
=P_{1}^{-+}P_{2}^{-+}...P_{t-1}^{-+}\left( P_{t}^{K}\right) ^{-+},\text{ }%
\forall i\in \left\langle n_{1}\right\rangle , 
\]
\noindent \textit{where} 
\[
\left( P_{t}^{K}\right) ^{-+}=\left( P_{t}^{K}\right) ^{-+\left( \Delta
_{t},\left( \left\{ j\right\} \right) _{j\in K}\right) } 
\]
\noindent ($\left( P_{1}P_{2}...P_{t}\right) _{\left\{ i\right\} }^{K}$ 
\textit{is the row }$i$\textit{\ of} \textit{matrix }$\left(
P_{1}P_{2}...P_{t}\right) ^{K}$; $\left(
P_{1}^{-+}P_{2}^{-+}...P_{t}^{-+}\right) ^{\bigcup\limits_{j\in K}\left\{
\left\{ j\right\} \right\} }$ \textit{is a submatrix of} $%
P_{1}^{-+}P_{2}^{-+}...P_{t}^{-+}$ (\textit{see the definition of operator} $%
\left( \cdot \right) ^{\left( \cdot \right) }$); $P_{b}^{-+}=P_{b}^{-+\left(
\Delta _{b},\Delta _{b+1}\right) },$ $\forall b\in \left\langle
t\right\rangle $; \textit{we used }$\left\{ \left\{ j\right\} \right\} $%
\textit{\ because the labels of columns of }$P_{1}^{-+}P_{2}^{-+}$

\noindent $...P_{t}^{-+}$\textit{\ are, here, }$\left\{ j\right\} ,$\textit{%
\ }$j\in \left\langle n_{t+1}\right\rangle ,$\textit{\ when }$\left|
K^{c}\right| \leq 1$ \textit{and }$\left\{ j\right\} ,$\textit{\ }$j\in K,$%
\textit{\ and }$K^{c}$ \textit{when }$\left| K^{c}\right| >1$\textit{\ }(%
\textit{see the definition of operator }$\left( \cdot \right) ^{-+}$ \textit{%
again}))\textit{, and} 
\[
\left( P_{1}P_{2}...P_{t}\right) ^{K}=e^{\prime }\left(
P_{1}^{-+}P_{2}^{-+}...P_{t}^{-+}\right) ^{\bigcup\limits_{j\in K}\left\{
\left\{ j\right\} \right\} }=e^{\prime
}P_{1}^{-+}P_{2}^{-+}...P_{t-1}^{-+}\left( P_{t}^{K}\right) ^{-+}, 
\]
\noindent \textit{where }$e=e\left( n_{1}\right) $ (\textit{if} $%
K=\left\langle n_{t+1}\right\rangle ,$\textit{\ then} 
\[
\left( P_{1}P_{2}...P_{t}\right) _{\left\{ i\right\}
}=P_{1}^{-+}P_{2}^{-+}...P_{t}^{-+},\text{ }\forall i\in \left\langle
n_{1}\right\rangle , 
\]
\noindent \textit{and} 
\[
P_{1}P_{2}...P_{t}=e^{\prime }P_{1}^{-+}P_{2}^{-+}...P_{t}^{-+}\text{).} 
\]

\smallskip

\textbf{Proof.} Since $\Delta _{1}=\left( \left\langle n_{1}\right\rangle
\right) $ and $\Delta _{t+1}=\left( \left\{ j\right\} ,K^{c}\right) _{j\in
K},$ by Theorem 1.9(ii) and Definition 1.5, $P_{1}P_{2}...P_{t}$ is a stable
matrix on $\left( \left\{ j\right\} ,K^{c}\right) _{j\in K}$ ($\left|
K^{c}\right| \leq 1$ or $\left| K^{c}\right| >1$). It follows that $\left(
P_{1}P_{2}...P_{t}\right) ^{K}$ is a stable matrix and, consequently, 
\[
\left( \left( P_{1}P_{2}...P_{t}\right) ^{K}\right) ^{-+\left( \left(
\left\langle n_{1}\right\rangle \right) ,\left( \left\{ j\right\} \right)
_{j\in K}\right) } 
\]
\noindent has only one row --- the label of this row is $\left\langle
n_{1}\right\rangle $ (see the definition of operator $\left( \cdot \right)
^{-+}$); further, we have 
\[
\left( P_{1}P_{2}...P_{t}\right) _{\left\{ i\right\} }^{K}=\left(
P_{1}P_{2}...P_{t-1}P_{t}^{K}\right) _{\left\{ i\right\} }= 
\]
\noindent (due to the partitions $\left( \left\langle n_{1}\right\rangle
\right) $ and $\left( \left\{ j\right\} \right) _{j\in K}...$ --- see
Example 1.8 again) 
\[
=\left( P_{1}P_{2}...P_{t-1}P_{t}^{K}\right) ^{-+\left( \left( \left\langle
n_{1}\right\rangle \right) ,\left( \left\{ j\right\} \right) _{j\in
K}\right) }= 
\]
\[
=P_{1}^{-+}P_{2}^{-+}...P_{t-1}^{-+}\left( P_{t}^{K}\right) ^{-+},\text{ }%
\forall i\in \left\langle n_{1}\right\rangle . 
\]
\noindent On the other hand, it is easy to see (see also Example 1.8) that 
\[
\left( P_{1}P_{2}...P_{t}\right) _{\left\{ i\right\} }^{K}=\left(
P_{1}^{-+}P_{2}^{-+}...P_{t}^{-+}\right) ^{\bigcup\limits_{j\in K}\left\{
\left\{ j\right\} \right\} },\text{ }\forall i\in \left\langle
n_{1}\right\rangle . 
\]

\noindent Finally, we have 
\[
\left( P_{1}P_{2}...P_{t}\right) ^{K}=e^{\prime }\left(
P_{1}^{-+}P_{2}^{-+}...P_{t}^{-+}\right) ^{\bigcup\limits_{j\in K}\left\{
\left\{ j\right\} \right\} }=e^{\prime
}P_{1}^{-+}P_{2}^{-+}...P_{t-1}^{-+}\left( P_{t}^{K}\right) ^{-+}.\text{ }%
\endproof%
\]

\smallskip

\textbf{Remark 1.11.} On Theorem 1.10, we consider 3 things.

\smallskip

(a) $\Delta _{t+1}=\left( \left\{ j\right\} \right) _{j\in \left\langle
n_{t+1}\right\rangle }$ both when $\left| K^{c}\right| =0$ and when $\left|
K^{c}\right| =1$ --- practically speaking, we use ``$\left| K^{c}\right| =0$%
'' because, in this case, $K=\left\langle n_{t+1}\right\rangle .$

\smallskip

(b) We have, in fact, 
\[
\left( P_{1}P_{2}...P_{t}\right) _{\left\{ i\right\} }^{K}\equiv \left(
P_{1}^{-+}P_{2}^{-+}...P_{t}^{-+}\right) ^{\bigcup\limits_{j\in K}\left\{
\left\{ j\right\} \right\} },\text{ }\forall i\in \left\langle
n_{1}\right\rangle , 
\]
\noindent because the row matrices 
\[
\left( P_{1}P_{2}...P_{t}\right) _{\left\{ i\right\} }^{K}\text{ and }\left(
P_{1}^{-+}P_{2}^{-+}...P_{t}^{-+}\right) ^{\bigcup\limits_{j\in K}\left\{
\left\{ j\right\} \right\} } 
\]
\noindent have different labels for rows and, separately, for columns, $%
\forall i\in \left\langle n_{1}\right\rangle ,$ but abusively, for
simplification, we used $``="$ instead of $``\equiv "$ ($\equiv $ is the
identical sign; since the one-to-one correspondence (the bijective function)
is column $j\longmapsto $ column (with label) $\left\{ j\right\} ,$ $\forall
j\in K,$ obviously, $\left( P_{1}^{-+}P_{2}^{-+}...P_{t}^{-+}\right)
^{\bigcup\limits_{j\in K}\left\{ \left\{ j\right\} \right\} }$ must have the
property: the column $\left\{ j_{1}\right\} $ is before the column $\left\{
j_{2}\right\} $ if $j_{1}<j_{2},$ $\forall j_{1},$ $j_{2}\in K$). Also, 
\[
\left( P_{1}P_{2}...P_{t}\right) _{\left\{ i\right\} }^{K}\equiv
P_{1}^{-+}P_{2}^{-+}...P_{t-1}^{-+}\left( P_{t}^{K}\right) ^{-+},\text{ }%
\forall i\in \left\langle n_{1}\right\rangle , 
\]
\noindent but, due to the notion of equal matrices and notation used for the
equal matrices, 
\[
\left( P_{1}^{-+}P_{2}^{-+}...P_{t}^{-+}\right) ^{\bigcup\limits_{j\in
K}\left\{ \left\{ j\right\} \right\}
}=P_{1}^{-+}P_{2}^{-+}...P_{t-1}^{-+}\left( P_{t}^{K}\right) ^{-+}. 
\]
\noindent Also, 
\[
\left( P_{1}P_{2}...P_{t}\right) ^{K}\equiv e^{\prime }\left(
P_{1}^{-+}P_{2}^{-+}...P_{t}^{-+}\right) ^{\bigcup\limits_{j\in K}\left\{
\left\{ j\right\} \right\} }... 
\]

\smallskip

(c) Theorem 1.10 leads to another one if we replace $\Delta _{t+1}=\left(
\left\{ j\right\} ,K^{c}\right) _{j\in K}$ with $\Delta _{t+1}=\left(
\left\{ j\right\} ,L_{2},L_{3},...,L_{v}\right) _{j\in L_{1}},$ where $%
\left( L_{1},L_{2},L_{3},...,L_{v}\right) \in $Par$\left( \left\langle
n_{t+1}\right\rangle \right) $ and $v\geq 2.$ But, by Remark 1.7(d), $%
P_{1}P_{2}...P_{t}$ is stable on $\left( \left\{ j\right\}
,L_{2},L_{3},...,L_{v}\right) _{j\in L_{1}}$

\noindent $\Longrightarrow P_{1}P_{2}...P_{t}$ is stable on$\left( \left\{
j\right\} ,\bigcup\limits_{b=2}^{v}L_{b}\right) _{j\in L_{1}}\!\!\!.$ So,
nothing new.

\smallskip

\textbf{Example 1.12.} Let 
\[
P_{1}=\left( 
\begin{array}{cccc}
\frac{2}{4}\smallskip & 0 & \frac{2}{4} & 0 \\ 
\frac{1}{4} & \frac{2}{4}\smallskip & \frac{1}{4} & 0 \\ 
0 & 0 & \frac{2}{4}\smallskip & \frac{2}{4} \\ 
0 & \frac{1}{4} & \frac{1}{4} & \frac{2}{4}%
\end{array}
\right) ,\text{ }P_{2}=\left( 
\begin{array}{cccc}
\frac{1}{4}\smallskip & \frac{1}{4} & 0 & \frac{2}{4} \\ 
\frac{2}{4} & 0\smallskip & 0 & \frac{2}{4} \\ 
\frac{2}{4} & 0 & \frac{1}{4}\smallskip & \frac{1}{4} \\ 
0 & \frac{2}{4} & \frac{2}{4} & 0%
\end{array}
\right) , 
\]
\[
P_{3}=\left( 
\begin{array}{cccc}
\frac{1}{4}\smallskip & 0 & \frac{2}{4} & \frac{1}{4} \\ 
\frac{1}{4} & 0\smallskip & \frac{2}{4} & \frac{1}{4} \\ 
0 & \frac{1}{4} & \frac{3}{4}\smallskip & 0 \\ 
0 & \frac{1}{4} & \frac{3}{4} & 0%
\end{array}
\right) ,\text{ and }P_{3}^{\prime }=\left( 
\begin{array}{cccc}
\frac{1}{4}\smallskip & 0 & \frac{2}{4} & \frac{1}{4} \\ 
\frac{1}{4} & 0\smallskip & \frac{1}{4} & \frac{2}{4} \\ 
0 & \frac{1}{4} & \frac{2}{4}\smallskip & \frac{1}{4} \\ 
0 & \frac{1}{4} & \frac{1}{4} & \frac{2}{4}%
\end{array}
\right) . 
\]
\noindent We have $P_{1}\in G_{\left( \left\langle 4\right\rangle \right)
,\left( \left\langle 4\right\rangle \right) },$ $P_{2}\in G_{\left(
\left\langle 4\right\rangle \right) ,\left( \left\{ 1,2\right\} ,\left\{
3,4\right\} \right) },$ $P_{3}\in G_{\left( \left\{ 1,2\right\} ,\left\{
3,4\right\} \right) ,\left( \left\{ j\right\} \right) _{j\in \left\langle
4\right\rangle }},$ and $P_{3}^{\prime }\in G_{\left( \left\{ 1,2\right\}
,\left\{ 3,4\right\} \right) ,\left( \left\{ 1\right\} ,\left\{ 2\right\}
,\left\{ 3,4\right\} \right) }.$ ($P_{1}\notin G_{\left( \left\langle
4\right\rangle \right) ,\Sigma },$ $\forall \Sigma \in $Par$\left(
\left\langle 4\right\rangle \right) ,$ $\Sigma \neq \left( \left\langle
4\right\rangle \right) $ --- $P_{1}$ was deliberately chosen with this
property.) So, by Theorem 1.10, $P_{1}P_{2}P_{3}$ is a stable matrix and 
\[
P_{1}P_{2}P_{3}=e^{\prime }P_{1}^{-+}P_{2}^{-+}P_{3}^{-+}=\left( 
\begin{array}{c}
1 \\ 
1 \\ 
1 \\ 
1%
\end{array}
\right) \left( 1\right) \left( 
\begin{array}{cc}
\frac{2}{4} & \frac{2}{4}%
\end{array}
\right) \left( 
\begin{array}{cccc}
\frac{1}{4}\smallskip & 0 & \frac{2}{4} & \frac{1}{4} \\ 
0 & \frac{1}{4} & \frac{3}{4} & 0%
\end{array}
\right) = 
\]
\[
=\left( 
\begin{array}{c}
1 \\ 
1 \\ 
1 \\ 
1%
\end{array}
\right) \left( 
\begin{array}{cccc}
\frac{2}{16} & \frac{2}{16} & \frac{10}{16} & \frac{2}{16}%
\end{array}
\right) =\left( 
\begin{array}{cccc}
\frac{1}{8}\smallskip & \frac{1}{8} & \frac{5}{8} & \frac{1}{8} \\ 
\frac{1}{8} & \frac{1}{8}\smallskip & \frac{5}{8} & \frac{1}{8} \\ 
\frac{1}{8} & \frac{1}{8} & \frac{5}{8}\smallskip & \frac{1}{8} \\ 
\frac{1}{8} & \frac{1}{8} & \frac{5}{8} & \frac{1}{8}%
\end{array}
\right) 
\]
\noindent while $\left( P_{1}P_{2}P_{3}^{\prime }\right) ^{\left\{
1,2\right\} }$ is a stable matrix and 
\[
\left( P_{1}P_{2}P_{3}^{\prime }\right) ^{\left\{ 1,2\right\} }=e^{\prime
}\left( P_{1}^{-+}P_{2}^{-+}P_{3}^{\prime -+}\right) ^{\left\{ \left\{
1\right\} ,\left\{ 2\right\} \right\} }=e^{\prime
}P_{1}^{-+}P_{2}^{-+}\left( P_{3}^{\prime -+}\right) ^{\left\{ \left\{
1\right\} ,\left\{ 2\right\} \right\} }= 
\]
\[
=\left( 
\begin{array}{c}
1 \\ 
1 \\ 
1 \\ 
1%
\end{array}
\right) \left( 1\right) \left( 
\begin{array}{cc}
\frac{2}{4} & \frac{2}{4}%
\end{array}
\right) \left( 
\begin{array}{cc}
\frac{1}{4}\smallskip & 0 \\ 
0 & \frac{1}{4}%
\end{array}
\right) =\left( 
\begin{array}{c}
1 \\ 
1 \\ 
1 \\ 
1%
\end{array}
\right) \left( 
\begin{array}{cc}
\frac{2}{16} & \frac{2}{16}%
\end{array}
\right) =\left( 
\begin{array}{cc}
\frac{1}{8}\smallskip & \frac{1}{8} \\ 
\frac{1}{8}\smallskip & \frac{1}{8} \\ 
\frac{1}{8}\smallskip & \frac{1}{8} \\ 
\frac{1}{8} & \frac{1}{8}%
\end{array}
\right) . 
\]
\noindent These results can also be obtained by direct computation. $%
P_{2}P_{3}$ and $\left( P_{2}P_{3}^{\prime }\right) ^{\left\{ 1,2\right\} }$
are also stable matrices.

\smallskip

To easy recognize when Theorem 1.10 could or cannot be used/applied, see (i)
and (ii) in the next theorem.

\smallskip

\textbf{THEOREM 1.13.} \textit{Suppose that the conditions of Theorem }1.10%
\textit{\ hold.}

\smallskip

(i) \textit{If }$n_{2},n_{3},...,n_{t+1}\geq 2$ ($n_{1}\geq 1$),\textit{\
then} $\exists u\in \left\langle t\right\rangle $\textit{\ such that }$P_{u}$%
\textit{\ is stable on }$\Delta _{u+1},$\textit{\ and }$\Delta _{u+1}\neq
\left( \left\langle n_{u+1}\right\rangle \right) $ ($P_{u}$ \textit{can be }%
``\textit{broken}''\textit{\ into at least two generalized stochastic
matrices with }$n_{u}$\textit{\ rows}).

\smallskip

(ii) \textit{If }$n_{1},n_{2},...,n_{t}\geq 2$ ($n_{t+1}\geq 1$),\textit{\
then} $\exists w\in \left\langle t-1\right\rangle $\textit{\ }($t\geq 2$) 
\textit{such that }$\Delta _{w}\neq \left( \left\{ i\right\} \right) _{i\in
\left\langle n_{w}\right\rangle }$ \textit{and} $\left( P_{w}\right) _{U}$%
\textit{\ is a matrix with identical rows, }$\forall U\in \Delta _{w},$ $%
\left| U\right| \geq 2$ (\textit{as a result}, $P_{w}$\textit{\ has at least
two identical rows})\textit{\ or }$\left( P_{t}\right) _{U}^{K}$\textit{\ is
a matrix with identical rows, }$\forall U\in \Delta _{t},$ $\left| U\right|
\geq 2$ (\textit{as a result}, $P_{t}^{K}$\textit{\ has at least two
identical rows }($P_{t}^{K}=P_{t}$\textit{\ when }$K=\left\langle
n_{t+1}\right\rangle $)).

\smallskip

(iii) $P_{1},$\textit{\ }$P_{1}P_{2},$\textit{\ }$...,$\textit{\ }$%
P_{1}P_{2}...P_{t}$\textit{\ are stable on }$\Delta _{2},$\textit{\ }$\Delta
_{3},$\textit{\ }$...,$\textit{\ }$\Delta _{t+1}$\textit{\ }($\Delta
_{t+1}=\left( \left\{ j\right\} ,K^{c}\right) _{j\in K}$), \textit{%
respectively, and, as a result, they are generalized stochastic matrices}.

\smallskip

(iv) $P_{t},$ $P_{t-1}P_{t},$ $...,$ $P_{1}P_{2}...P_{t}$ \textit{are }$%
\left[ \Delta _{t}\right] $\textit{-, }$\left[ \Delta _{t-1}\right] $\textit{%
-, ..., }$\left[ \Delta _{1}\right] $\textit{-stable on }$\left( \left\{
j\right\} ,K^{c}\right) _{j\in K}$,\textit{\ respectively} ($\Delta
_{1}=\left( \left\langle n_{1}\right\rangle \right) $).

\smallskip

\textbf{Proof.} (i) Since $n_{2},n_{3},...,n_{t+1}\geq 2,$\textit{\ }$\Delta
_{1}\!=\!\left( \left\langle n_{1}\right\rangle \right) ,$ and $\Delta
_{t+1}\!=\!\left( \left\{ j\right\} ,K^{c}\right) _{j\in K}\!,$ it follows
that $\exists u\in \left\langle t\right\rangle $ such that $\Delta
_{1}=\left( \left\langle n_{1}\right\rangle \right) $, $\Delta _{2}=\left(
\left\langle n_{2}\right\rangle \right) $, ..., $\Delta _{u}=\left(
\left\langle n_{u}\right\rangle \right) $, and $\Delta _{u+1}\neq \left(
\left\langle n_{u+1}\right\rangle \right) $ ($``u=t"$ leads to $``\Delta
_{t+1}\neq \left( \left\langle n_{t+1}\right\rangle \right) "$; $n_{t+1}\geq
2$ and $\emptyset \neq K\subseteq \left\langle n_{t+1}\right\rangle
\Longrightarrow \left| \Delta _{t+1}\right| \geq 2$ (equivalently, $\Delta
_{t+1}\neq \left( \left\langle n_{t+1}\right\rangle \right) $); so, no
problem). It follows that $P_{u}$ is stable on $\Delta _{u+1},$ and $\Delta
_{u+1}\neq \left( \left\langle n_{u+1}\right\rangle \right) $.

\smallskip

(ii) \textit{Case} 1. $\Delta _{t}\neq \left( \left\{ i\right\} \right)
_{i\in \left\langle n_{t}\right\rangle }$ --- this case is possible because $%
n_{t}\geq 2.$ In this case, $\exists V\in \Delta _{t}$ with $\left| V\right|
\geq 2.$ Since $P_{t}\in \widehat{G}_{\Delta _{t},\Delta _{t+1}}$ and $%
\Delta _{t+1}=\left( \left\{ j\right\} ,K^{c}\right) _{j\in K},$ it follows
that $\left( P_{t}\right) _{U}^{K}$ is a (complex) matrix with identical
rows, $\forall U\in \Delta _{t},$ $\left| U\right| \geq 2.$

\smallskip

\textit{Case} 2. $\Delta _{t}=\left( \left\{ i\right\} \right) _{i\in
\left\langle n_{t}\right\rangle }.$ Since $n_{1},n_{2},...,n_{t-1}\geq 2$
(here, $n_{t-1},$ not $n_{t}$), $\Delta _{1}=\left( \left\langle
n_{1}\right\rangle \right) ,$ and $\Delta _{t}=\left( \left\{ i\right\}
\right) _{i\in \left\langle n_{t}\right\rangle },$ it follows that $\exists
w\in \left\langle t-1\right\rangle $ ($t\geq 2$) such that $\Delta
_{t}=\left( \left\{ i\right\} \right) _{i\in \left\langle n_{t}\right\rangle
},$ $\Delta _{t-1}=\left( \left\{ i\right\} \right) _{i\in \left\langle
n_{t-1}\right\rangle },$ ..., $\Delta _{w+1}=\left( \left\{ i\right\}
\right) _{i\in \left\langle n_{w+1}\right\rangle },$ and $\Delta _{w}\neq
\left( \left\{ i\right\} \right) _{i\in \left\langle n_{w}\right\rangle }$ ($%
n_{w}\geq 2,$ so, $``\Delta _{w}\neq \left( \left\{ i\right\} \right) _{i\in
\left\langle n_{w}\right\rangle }"$ is possible). Since $P_{w}\in \widehat{G}%
_{\Delta _{w},\Delta _{w+1}},$ it follows that $\left( P_{w}\right) _{U}$ is
a matrix with identical rows, $\forall U\in \Delta _{w},$ $\left| U\right|
\geq 2$ ($\Delta _{w}\neq \left( \left\{ i\right\} \right) _{i\in
\left\langle n_{w}\right\rangle }\Longrightarrow \exists V\in \Delta _{w}$
with $\left| V\right| \geq 2$).

\smallskip

(iii) and (iv) Theorem 1.9. 
\endproof%

\smallskip

Theorem 1.10 was illustrated in Example 1.12. To illustrate Theorem
1.13(i)-(ii), considering Example 1.12, the matrix $P_{1}$ cannot be
``broken'' into at least two generalized stochastic matrices with 4 rows,
the matrix $P_{2}$ can be ``broken'' into two generalized stochastic
(nonnegative) matrices with 4 rows, namely, $P_{2}^{\left\langle
2\right\rangle }$ and $P_{2}^{\left\{ 3,4\right\} },$ 
\[
P_{2}^{\left\langle 2\right\rangle }=\left( 
\begin{array}{cc}
\frac{1}{4}\smallskip & \frac{1}{4} \\ 
\frac{2}{4}\smallskip & 0 \\ 
\frac{2}{4}\smallskip & 0 \\ 
0 & \frac{2}{4}%
\end{array}
\right) \text{ and }P_{2}^{\left\{ 3,4\right\} }=\left( 
\begin{array}{cc}
0\smallskip & \frac{2}{4} \\ 
0\smallskip & \frac{2}{4} \\ 
\frac{1}{4}\smallskip & \frac{1}{4} \\ 
\frac{2}{4} & 0%
\end{array}
\right) , 
\]
\noindent the matrix $P_{3}$ is $\left[ \left( \left\langle 2\right\rangle
,\left\{ 3,4\right\} \right) \right] $-stable, so, $\left( P_{3}\right)
_{\left\langle 2\right\rangle }$ and $\left( P_{3}\right) _{\left\{
3,4\right\} }$ are stable matrices 
\[
\text{(}\left( P_{3}\right) _{\left\langle 2\right\rangle }=\left( 
\begin{array}{cccc}
\frac{1}{4}\smallskip & 0 & \frac{2}{4} & \frac{1}{4} \\ 
\frac{1}{4} & 0 & \frac{2}{4} & \frac{1}{4}%
\end{array}
\right) ,...\text{),} 
\]
\noindent and the matrix $P_{3}^{\left\langle 2\right\rangle }$ is $\left[
\left( \left\langle 2\right\rangle ,\left\{ 3,4\right\} \right) \right] $%
-stable, so, $\left( P_{3}\right) _{\left\langle 2\right\rangle
}^{\left\langle 2\right\rangle }$ and $\left( P_{3}\right) _{\left\{
3,4\right\} }^{\left\langle 2\right\rangle }$ are stable matrices.

\smallskip

Below we consider a similarity relation, and give a result about it --- the
relation and result lead, under certain conditions, to a subset/subgroup
property of the DeGroot submodel from Section 2, see Example 2.11.

\smallskip

\textbf{Definition 1.14.} (A generalization of Definition 2.1 from [19].)
Let $P,Q\in \widehat{G}_{\Delta ,\Sigma }\subseteq C_{m,n}.$ We say that $P$%
\textit{\ is similar to }$Q$ (\textit{with respect to }$\Delta $\textit{\
and }$\Sigma $) if 
\[
P^{-+}=Q^{-+}. 
\]

Set $P\backsim Q$ when $P$ is similar to $Q.$ Obviously, $\backsim $ is an
equivalence relation on $\widehat{G}_{\Delta ,\Sigma }.$

\smallskip

\textbf{Example 1.15.} (For another example, see Example 1.1 in [21].) Let 
\[
P=\left( 
\begin{array}{cccc}
\frac{1}{4}\smallskip & 0 & \frac{2}{4} & \frac{1}{4} \\ 
0 & \frac{1}{4}\smallskip & \frac{3}{4} & 0 \\ 
\frac{1}{4} & \frac{1}{4} & \frac{1}{4}\smallskip & \frac{1}{4} \\ 
\frac{1}{4} & \frac{1}{4} & \frac{1}{4} & \frac{1}{4}\smallskip%
\end{array}
\right) ,\text{ }Q=\left( 
\begin{array}{cccc}
\frac{1}{4}\smallskip & 0 & \frac{3}{4} & 0 \\ 
\frac{1}{4} & 0\smallskip & 0 & \frac{3}{4} \\ 
\frac{2}{4} & 0 & \frac{2}{4}\smallskip & 0 \\ 
0 & \frac{2}{4} & \frac{2}{4} & 0\smallskip%
\end{array}
\right) ,\text{ and }T=\left( 
\begin{array}{cccc}
\frac{1}{4}\smallskip & 0 & \frac{3}{4} & 0 \\ 
\frac{1}{4} & 0\smallskip & \frac{2}{4} & \frac{1}{4} \\ 
0 & \frac{2}{4} & \frac{2}{4}\smallskip & 0 \\ 
\frac{2}{4} & 0 & \frac{1}{4} & \frac{1}{4}\smallskip%
\end{array}
\right) . 
\]
\noindent We have $P,Q,R\in G_{\Delta ,\Sigma }\subseteq S_{4},$ where $%
\Delta =\Sigma =\left( \left\{ 1,2\right\} ,\left\{ 3,4\right\} \right) ,$
and 
\[
P^{-+}=Q^{-+}=T^{-+}=\left( 
\begin{array}{cc}
\frac{1}{4}\smallskip & \frac{3}{4} \\ 
\frac{2}{4} & \frac{2}{4}%
\end{array}
\right) . 
\]
\noindent So, $P\backsim Q\backsim T.$ Note that: 1) any matrix which is
similar to $P$ with respect to $\Delta $ and $\Sigma $ has at least $8$
positive entries; 3) $Q$ has the smallest possible number of positive
entries, \textit{i.e.}, $8,$ while $T$ has not --- the matrices with a small
number of positive entries are important for algorithms (for sampling...).

\smallskip

The next result says, among other things, that under certain conditions a
product of matrices can be replaced with another product of matrices.

\smallskip

\textbf{THEOREM\ 1.16.} (A generalization of Theorem 2.3 from [19].) \textit{%
Let }$P_{1},U_{1}\in \widehat{G}_{\Delta _{1},\Delta _{2}}\subseteq
C_{n_{1},n_{2}},$\textit{\ }$P_{2},U_{2}\in \widehat{G}_{\Delta _{2},\Delta
_{3}}\subseteq C_{n_{2},n_{3}},$\textit{\ }$...,$\textit{\ }$P_{t},U_{t}\in 
\widehat{G}_{\Delta _{t},\Delta _{t+1}}\subseteq C_{n_{t},n_{t+1}}.$\textit{%
\ Suppose that} 
\[
P_{1}\backsim U_{1},\mathit{\ }P_{2}\backsim U_{2},\mathit{\ }%
...,P_{t}\backsim U_{t}. 
\]
\noindent \textit{Then} 
\[
P_{1}P_{2}...P_{t}\backsim U_{1}U_{2}....U_{t}. 
\]

\noindent \textit{If, moreover, }$\Delta _{1}=\left( \left\langle
n_{1}\right\rangle \right) $\textit{\ and }$\Delta _{t+1}=\left( \left\{
j\right\} ,K^{c}\right) _{j\in K},$\textit{\ where }$\emptyset \neq
K\subseteq \left\langle n_{t+1}\right\rangle $\textit{, then} 
\[
\left( P_{1}P_{2}...P_{t}\right) ^{K}=\left( U_{1}U_{2}....U_{t}\right) ^{K} 
\]
\[
\text{--- \textit{equivalently}, }%
P_{1}P_{2}...P_{t-1}P_{t}^{K}=U_{1}U_{2}....U_{t-1}U_{t}^{K} 
\]

\noindent (\textit{therefore, when }$\Delta _{1}=\left( \left\langle
n_{1}\right\rangle \right) $\textit{\ and }$\Delta _{t+1}=\left( \left\{
j\right\} ,K^{c}\right) _{j\in K},$\textit{\ a product of }$n$\textit{\
representatives, the first of an equivalence class included in }$\widehat{G}%
_{\Delta _{1},\Delta _{2}},$\textit{\ the second of an equivalence class
included in }$\widehat{G}_{\Delta _{2},\Delta _{3}},$\textit{\ }$...,$%
\textit{\ the }$t$th\textit{\ of an equivalence class included in }$\widehat{%
G}_{\Delta _{t},\left( \left\{ j\right\} \right) _{j\in K}},$\textit{\ does
not depend on the choice of representatives --- we can work }(\textit{this
is important})\textit{\ with }$U_{1}U_{2}....U_{t-1}U_{t}^{K}$ \textit{%
instead of} $P_{1}P_{2}...P_{t-1}P_{t}^{K}$ \textit{and}, \textit{%
conversely, with} $P_{1}P_{2}...P_{t-1}P_{t}^{K}$ \textit{instead of} $%
U_{1}U_{2}....U_{t-1}U_{t}^{K}$; \textit{due to Theorem} 1.10, $%
P_{1}P_{2}...P_{t-1}P_{t}^{K}$ \textit{and} $U_{1}U_{2}....U_{t-1}U_{t}^{K}$ 
\textit{are stable matrices};\textit{\ }$P_{1}P_{2}...P_{t-1}$

\noindent $P_{t}^{K}=P_{1}P_{2}...P_{t}$ \textit{and} $%
U_{1}U_{2}....U_{t-1}U_{t}^{K}=U_{1}U_{2}....U_{t}$ \textit{when} $%
K=\left\langle n_{t+1}\right\rangle $).

\smallskip

\textbf{Proof.} The first part is similar to the proof of first part of
Theorem 2.3 from [19]. For the second part, we consider two cases.

\textit{Case} 1. $\left| K^{c}\right| \leq 1.$ The proof is similar to the
proof of second part of Theorem 2.3 from [19].

\textit{Case} 2. $\left| K^{c}\right| >1.$ We have (see Remark 1.2) 
\[
\left( P_{1}P_{2}...P_{t}\right) ^{K}=P_{1}P_{2}...P_{t-1}P_{t}^{K} 
\]
\noindent and 
\[
\left( U_{1}U_{2}....U_{t}\right) ^{K}=U_{1}U_{2}....U_{t-1}U_{t}^{K}. 
\]
\noindent Since $P_{t},U_{t}\in \widehat{G}_{\Delta _{t},\left( \left\{
j\right\} ,K^{c}\right) _{j\in K}},$ we have $P_{t}^{K},U_{t}^{K}\in 
\widehat{G}_{\Delta _{t},\left( \left\{ j\right\} \right) _{j\in K}}.$ To
finish the proof, we apply Case 1 to the matrices $P_{1},$ $P_{2},$ $...,$ $%
P_{t-1},$ $P_{t}^{K}$ and $U_{1},$ $U_{2},$ $....,$ $U_{t-1},$ $U_{t}^{K},$
and obtain 
\[
P_{1}P_{2}...P_{t-1}P_{t}^{K}=U_{1}U_{2}....U_{t-1}U_{t}^{K}.\text{ }%
\endproof%
\]

\smallskip

For examples for Theorem 1.16, see Example 1.2 in [21] and, in this article,
Example 2.11 (in Section 2).

\smallskip

\textbf{Remark 1.17.} The equivalence relation $\backsim $ is, in
particular, an equivalence relation on the set of generalized stochastic
complex $m\times n$ matrices. Indeed, by Remark 1.7(b) we have that $%
\widehat{G}_{\left( \left\langle m\right\rangle \right) ,\left( \left\langle
n\right\rangle \right) }$ is equal to the set of generalized stochastic
complex $m\times n$ matrices. Since $\backsim $ is an equivalence relation
(with respect to $\left( \left\langle m\right\rangle \right) $ and $\left(
\left\langle n\right\rangle \right) $) on $\widehat{G}_{\left( \left\langle
m\right\rangle \right) ,\left( \left\langle n\right\rangle \right) }$, it
follows that it is an equivalence relation (with respect to $\left(
\left\langle m\right\rangle \right) $ and $\left( \left\langle
n\right\rangle \right) $ --- see Definition 1.14 again) on the set of
generalized stochastic complex $m\times n$ matrices. Similarly, $\backsim $
is an equivalence relation on the set of generalized stochastic real $%
m\times n$ matrices and on that of generalized stochastic nonnegative $%
m\times n$ matrices.

\bigskip

\begin{center}
\textbf{2. FINITE-TIME\ CONSENSUS\ FOR\ DEGROOT\ MODEL}
\end{center}

\bigskip

In this section, we consider the (homogeneous and nonhomogeneous) DeGroot
model and the finite-time consensus problem for this model --- using the $G$
method, a result for reaching a partial or total consensus in finite time is
given. Then we consider a special submodel/case of the DeGroot model, and
for it we give two examples, and for the latter example (Example 2.11) we
give certain comments --- a subset/subgroup property of the DeGroot submodel
is illustrated in this latter example.

\smallskip

Consider a group (society, team, committee, etc.) of $r$ individuals
(experts, agents, etc.) --- formally/mathematically we consider a set with $%
r $ elements, $r\geq 2$; in this article, we consider the set $\left\langle
r\right\rangle $ ($\left\langle r\right\rangle =\left\{ 1,2,...,r\right\} ,$ 
$r\geq 2$), and we call it the\textit{\ individual set}. Each individual has
an estimation/opinion, an initial one, on a subject (the estimation of a
parameter, or of a probability, or...). The initial estimations/opinions of
individuals are represented by a (row) vector $p_{0}=\left( p_{0i}\right)
_{i\in \left\langle r\right\rangle }=\left( p_{01},p_{02},...,p_{0r}\right) $%
, $p_{0i}$ is the initial estimation/opinion of individual $i,$ $\forall
i\in \left\langle r\right\rangle ,$ and $p_{0}$ is a probability (row)
vector, or a nonnegative vector, or a real vector, or a complex vector,
or... When the individuals are made aware of each others' opinions, they may
modify their own opinion by taking into account the opinion of others, and
the modifications, if any, are done at discrete times $1,2,...$ ---
formally, we consider a sequence $\left( p_{n}\right) _{n\geq 0}$ of (row)
vectors and a sequence $\left( P_{n}\right) _{n\geq 1}$ of stochastic
matrices, $p_{n}=\left( p_{ni}\right) _{i\in \left\langle r\right\rangle
}=\left( p_{n1},p_{n2},...,p_{nr}\right) ,$ $\forall n\geq 0,$ and $%
P_{n}=\left( \left( P_{n}\right) _{ij}\right) \in S_{r},$ $\forall n\geq 1,$
where $p_{ni}$ is the estimation/opinion of individual $i$ at time $n,$ $%
\forall n\geq 0,$ $\forall i\in \left\langle r\right\rangle ,$ and $\left(
P_{n}\right) _{ij}$ is the probability which the individual $i$ assigns to
the estimation/opinion $p_{n-1j}$ (warning! not $p_{0j}$) of individual $j$
at time $n$ (not $n-1$), $\forall n\geq 1,$ $\forall i,j\in \left\langle
r\right\rangle .$ We call $p_{n}$ the \textit{vector of estimations/opinions 
}(\textit{of individuals})\textit{\ at time }$n$\textit{\ or the value
vector at time }$n,$ $\forall n\geq 0,$ and $P_{n}$ the \textit{trust }or%
\textit{\ weight matrix} \textit{at time} $n,$ $\forall n\geq 1.$ We also
call $p_{0}$ the\textit{\ initial vector of estimations/opinions }(\textit{%
of individuals})\textit{\ or the initial value vector. }$p_{n}$ is a
probability (row) vector, or a nonnegative vector, or a real vector, or a
complex vector, or..., $\forall n\geq 0.$ The above considerations lead to 
\[
p_{n}^{\prime }=P_{n}p_{n-1}^{\prime },\text{ }\forall n\geq 1 
\]
\noindent ($p_{n}^{\prime }$ and $p_{n-1}^{\prime }$ are the transposes of $%
p_{n}$ and $p_{n-1}$, respectively, $\forall n\geq 1$), and, as a result, to 
\[
p_{n}^{\prime }=P_{n}P_{n-1}...P_{1}p_{0}^{\prime },\text{ }\forall n\geq 1. 
\]

\smallskip

The above model is determined by the set $\left\langle r\right\rangle ,$
vector $p_{0},$ and stochastic matrices $P_{n},$ $n\geq 1.$ Due this fact,
we call this model the \textit{DeGroot model on }(\textit{the triple}) $%
\left( \left\langle r\right\rangle ,p_{0},\left( P_{n}\right) _{n\geq
1}\right) ,$ and will use the generic name ``(the) DeGroot model''. We say
that DeGroot model on $\left( \left\langle r\right\rangle ,p_{0},\left(
P_{n}\right) _{n\geq 1}\right) $ is \textit{homogeneous} if $%
P_{1}=P_{2}=P_{3}=...$ and \textit{nonhomogeneous} if $\exists u,v\geq 1,$ $%
u\neq v,$ such that $P_{u}\neq P_{v},$ and will use the generic names
``(the) homogeneous DeGroot model'' and ``(the) nonhomogeneous DeGroot
model'', respectively.

\smallskip

\textbf{Definition 2.1} ([7]; see, \textit{e.g.}, also [5]). We say that the 
\textit{individual set} \textit{reaches a consensus} or that a \textit{%
consensus is reached }(\textit{for the individual set}) if all $r$
components of $p_{n}$ converge to the same limit as $n\longrightarrow \infty 
$, and we call this limit the \textit{consensus}.

\smallskip

Reaching a consensus for the DeGroot model --- this is the \textit{consensus
problem} for/of the DeGroot model. When a consensus is reached, it is
reached in a finite or infinite time --- obviously, the first case is the
best.

\smallskip

The DeGroot model and its consensus problem and the consensus problem in
computer science have some things in common with the finite Markov chains
--- for the first three, see, \textit{e.g.}, [1, Chapter 5, pp. 91$-$124],
[2]-[3], [5]-[7], [9], [12, Chapter 8], [15], [22], [26]-[28], and [30], and
for the last, see, \textit{e.g.}, [10]-[11], [13], [16], and [24]-[26] (all
or some of these references could be available; arXiv and Wikipedia could
also be useful... --- expressions for search: DeGroot model, DeGroot
learning, DeGroot learning model, consensus (computer science), distributed
consensus, distributed averaging, Markov chain, homogeneous Markov chain,
inhomogeneous Markov chain, nonhomogeneous Markov chain, ...).

\smallskip

\textbf{THEOREM\ 2.2.} (A simple, important, and known result, see [7], [2],
...) \textit{Consider the DeGroot model on }$\left( \left\langle
r\right\rangle ,p_{0},\left( P_{n}\right) _{n\geq 1}\right) .$\textit{\ If }$%
\lim\limits_{n\rightarrow \infty }P_{n}P_{n-1}...P_{1}$\textit{\ exists and
is a stable matrix, then a consensus is reached.}

\smallskip

\textbf{Proof.} (A known proof, but we give it --- the best thing is to give
it.) Suppose that $\lim\limits_{n\rightarrow \infty }P_{n}P_{n-1}...P_{1}$
exists and is a stable matrix. In this case, $\exists \pi ,$ $\pi $ is a
probability (row) vector with $r$ components, such that 
\[
\lim\limits_{n\rightarrow \infty }P_{n}P_{n-1}...P_{1}=e^{\prime }\pi 
\]
\noindent (recall that $e=e\left( r\right) =\left( 1,1,...,1\right) $).
Further, we have 
\[
\lim\limits_{n\rightarrow \infty }p_{n}^{\prime }=\lim\limits_{n\rightarrow
\infty }P_{n}P_{n-1}...P_{1}p_{0}^{\prime }=e^{\prime }\pi p_{0}^{\prime
}=e^{\prime }\sum\limits_{i=1}^{r}\pi _{i}p_{0i}=\left( 
\begin{array}{c}
\sum\limits_{i=1}^{r}\pi _{i}p_{0i} \\ 
\sum\limits_{i=1}^{r}\pi _{i}p_{0i} \\ 
\vdots \\ 
\sum\limits_{i=1}^{r}\pi _{i}p_{0i}%
\end{array}
\right) . 
\]
\noindent Therefore, a consensus is reached; the consensus is 
\[
\sum\limits_{i=1}^{r}\pi _{i}p_{0i}.\text{ }%
\endproof%
\]

\smallskip

\textbf{Remark 2.3.} The condition from the above theorem is only a
sufficient condition for reaching a consensus. Indeed, if 
\[
p_{0}=\left( 2,4,3\right) \text{ and }P_{n}=\left( 
\begin{array}{ccc}
\frac{1}{2}\smallskip & \frac{1}{2} & 0 \\ 
\frac{1}{2} & \frac{1}{2}\smallskip & 0 \\ 
0 & 0 & 1%
\end{array}
\right) ,\forall n\geq 1 
\]
\noindent (a homogeneous case), then a consensus is reached at time $1$ (a
finite-time consensus), the consensus is $3$, but 
\[
\lim_{n\rightarrow \infty }P_{n}P_{n-1}...P_{1}=\left( 
\begin{array}{ccc}
\frac{1}{2}\smallskip & \frac{1}{2} & 0 \\ 
\frac{1}{2} & \frac{1}{2}\smallskip & 0 \\ 
0 & 0 & 1%
\end{array}
\right) , 
\]
\noindent which is not a stable matrix. Another example, a nonhomogeneous
case: 
\[
p_{0}=\left( 2,4,\frac{10}{3}\right) \text{, }P_{1}=\left( 
\begin{array}{ccc}
\frac{1}{3}\smallskip & \frac{2}{3} & 0 \\ 
\frac{1}{3} & \frac{2}{3}\smallskip & 0 \\ 
0 & 0 & 1%
\end{array}
\right) ,\text{ and }P_{n}=\left( 
\begin{array}{cc}
Q_{n} & 
\begin{array}{c}
0 \\ 
0%
\end{array}
\\ 
\begin{array}{cc}
0 & 0%
\end{array}
& 1%
\end{array}
\right) ,\forall n\geq 2, 
\]
\noindent where $Q_{n}$ is a stochastic $2\times 2$ matrix, $\forall n\geq 2$
--- a consensus is reached at time $1$, and this is $\frac{10}{3}.$ Note
that in the above cases/examples the matrices are reducible --- therefore, a
consensus, even a finite-time consensus, can be reached even the matrices
are reducible.

\smallskip

The above remark is both for Theorem 2.2 and for an incorrect statement in
[7, Sec. 4, p. 119] discovered by R.L. Berger [2].

\smallskip

\textbf{Remark 2.4.} (a) If $P\in S_{m,n}$ and $Q\in S_{n,p}$, and $P$ is a
stable matrix, then $PQ$ is a stable matrix. More generally, if $P\in
C_{m,n} $ and $Q\in C_{n,p}$, and $P$ is a stable (complex) matrix, then $PQ$
is a stable matrix.

\smallskip

(b) If $P\in S_{m,n}$ and $Q\in C_{n,p}$, and $Q$ is a stable matrix, then $%
PQ$ is a stable matrix, and, moreover, 
\[
PQ=Q. 
\]
\noindent More generally, if $P\in C_{m,n}$ and $Q\in C_{n,p}$, $P$ is a
generalized stochastic complex matrix (see Definition 1.3) and $Q$ is a
stable (complex) matrix, then $PQ$ is a stable matrix, and, moreover, 
\[
PQ=Q, 
\]
\noindent but only when $P$ is a stochastic complex matrix.

\smallskip

(c) (simple and important proprieties of the (homogeneous and
nonhomogeneous) DeGroot model) Consider the DeGroot model on $\left(
\left\langle r\right\rangle ,p_{0},\left( P_{n}\right) _{n\geq 1}\right) .$

\smallskip

\hspace{1cm}(c1) If $P_{t}P_{t-1}...P_{m}$ is a stable matrix for some $t$
and $m,$ $t\geq m\geq 1,$ then (see (a) and (b)) 
\[
P_{n}P_{n-1}...P_{1} 
\]
\noindent is a stable matrix, $\forall n\geq t,$ and, moreover, 
\[
P_{n}P_{n-1}...P_{1}=P_{t}P_{t-1}...P_{1},\text{ }\forall n\geq t. 
\]

\noindent This fact leads to other ones: 1) the sequence of matrices $\left(
P_{n}\right) _{n\geq t+1}$ does not count, it can be replaced with any other
infinite sequence of stochastic $r\times r$ matrices; 2) a consensus is
reached in a finite time, and it is 
\[
\sum\limits_{i=1}^{r}\pi _{i}p_{0i}, 
\]
\noindent where $\pi =\left( \pi _{1},\pi _{2},...,\pi _{r}\right) ,$ $\pi
=\left( P_{t}P_{t-1}...P_{1}\right) _{\left\{ i\right\} },$ $\forall i\in
\left\langle r\right\rangle $ ($P_{t}P_{t-1}...P_{1}$ is a stable matrix, 
\textit{i.e.}, a matrix with identical rows; $\left(
P_{t}P_{t-1}...P_{1}\right) _{\left\{ i\right\} }$ is the row $i$ of $%
P_{t}P_{t-1}...P_{1},$ $\forall i\in \left\langle r\right\rangle $).

\smallskip

\hspace{1cm}(c2) If $\left( P_{t}P_{t-1}...P_{1}\right) ^{K}$ is a stable
matrix for some $t$ and $K,$ $t\geq 1$ and $\emptyset \neq K\subseteq
\left\langle r\right\rangle ,$ then (see Remark 1.2(b) and (b)) 
\[
\left( P_{n}P_{n-1}...P_{1}\right) ^{K} 
\]
\noindent is a stable matrix, $\forall n\geq t,$ and, moreover, 
\[
\left( P_{n}P_{n-1}...P_{1}\right) ^{K}=\left( P_{t}P_{t-1}...P_{1}\right)
^{K},\text{ }\forall n\geq t 
\]
\noindent --- $\forall n\geq t,$ $P_{n}P_{n-1}...P_{1}=P_{t}P_{t-1}...P_{1}$
and $P_{n}P_{n-1}...P_{1}$ is a stable matrix when $\left| K^{c}\right| =0;$ 
$\forall n\geq t,$ $\left( P_{n}P_{n-1}...P_{1}\right) ^{K}=\left(
P_{t}P_{t-1}...P_{1}\right) ^{K}$ and $\left( P_{n}P_{n-1}...P_{1}\right)
^{K}$ is a stable matrix $\Longrightarrow
P_{n}P_{n-1}...P_{1}=P_{t}P_{t-1}...P_{1}$ and $P_{n}P_{n-1}...P_{1}$ is a
stable matrix when $\left| K^{c}\right| =1$ because the matrices $P_{b},$ $%
b\geq 1,$ are stochastic; so, $\forall n\geq t,$ $%
P_{n}P_{n-1}...P_{1}=P_{t}P_{t-1}...P_{1}$ and $P_{n}P_{n-1}...P_{1}$ is a
stable matrix when $\left| K^{c}\right| \leq 1.$

\smallskip

Below we give a fundamental result concerning the finite-time partial and
total consensuses (for reached a partial or total consensus in a finite
time).

\smallskip

\textbf{THEOREM 2.5.} \textit{Consider the DeGroot model} \textit{on }$%
\left( \left\langle r\right\rangle ,p_{0},\left( P_{n}\right) _{n\geq
1}\right) .$ \textit{If} $P_{1}\in G_{\Delta _{2},\Delta _{1}}\subseteq
S_{r},$ $P_{2}\in G_{\Delta _{3},\Delta _{2}}\subseteq S_{r}$ $...,$ $%
P_{t}\in G_{\Delta _{t+1},\Delta _{t}}\subseteq S_{r}$, \textit{and }$\Delta
_{t+1}=\left( \left\langle r\right\rangle \right) $ \textit{for some }$t\geq
1$\textit{\ and} 
\[
\Delta _{1}=\left( \left\{ j\right\} ,K^{c}\right) _{j\in K}\text{ }%
(\emptyset \neq K\subseteq \left\langle r\right\rangle ), 
\]
\noindent \textit{then} 
\[
\left( P_{t}P_{t-1}...P_{1}\right) ^{K} 
\]
\noindent (\textit{equivalently}, $P_{t}P_{t-1}...P_{2}P_{1}^{K}$) \textit{%
is a stable matrix}, 
\[
\left( P_{t}P_{t-1}...P_{1}\right) _{\left\{ i\right\} }^{K}\!=\!\left(
P_{t}^{-+}P_{t-1}^{-+}...P_{1}^{-+}\right) ^{\bigcup\limits_{j\in K}\left\{
\left\{ j\right\} \right\} }\!=\!P_{t}^{-+}P_{t-1}^{-+}...P_{2}^{-+}\left(
P_{1}^{K}\right) ^{-+}\!,\forall i\!\in \!\left\langle r\right\rangle \!\!, 
\]
\noindent \textit{where} 
\[
\left( P_{1}^{K}\right) ^{-+}=\left( P_{1}^{K}\right) ^{-+\left( \Delta
_{2},\left( \left\{ j\right\} \right) _{j\in K}\right) } 
\]
\noindent ($\left( P_{t}P_{t-1}...P_{1}\right) _{\left\{ i\right\} }^{K}$ 
\textit{is the row }$i$\textit{\ of} \textit{matrix }$\left(
P_{t}P_{t-1}...P_{1}\right) ^{K}$; $\left(
P_{t}^{-+}P_{t-1}^{-+}...P_{1}^{-+}\right) ^{\bigcup\limits_{j\in K}\left\{
\left\{ j\right\} \right\} }$ \textit{is a submatrix of} $%
P_{t}^{-+}P_{t-1}^{-+}...P_{1}^{-+}$ (\textit{see the definition of operator 
}$\left( \cdot \right) ^{\left( \cdot \right) }$\textit{\ again}); $%
P_{b}^{-+}=P_{b}^{-+\left( \Delta _{b+1},\Delta _{b}\right) },$ $\forall
b\in \left\langle t\right\rangle ;$ \textit{we used }$\left\{ \left\{
j\right\} \right\} $\textit{\ because the labels of columns of }$%
P_{t}^{-+}P_{t-1}^{-+}...P_{1}^{-+}$\textit{\ are, here, }$\left\{ j\right\}
,$\textit{\ }$j\in \left\langle r\right\rangle ,$\textit{\ when }$\left|
K^{c}\right| \leq 1$ \textit{and }$\left\{ j\right\} ,$\textit{\ }$j\in K,$%
\textit{\ and }$K^{c}$ \textit{when }$\left| K^{c}\right| >1$\textit{\ }(%
\textit{see the definition of operator }$\left( \cdot \right) ^{-+}$\textit{%
\ again}))\textit{, and} 
\[
\left( P_{t}P_{t-1}...P_{1}\right) ^{K}=e^{\prime }\left(
P_{t}^{-+}P_{t-1}^{-+}...P_{1}^{-+}\right) ^{\bigcup\limits_{j\in K}\left\{
\left\{ j\right\} \right\} }= 
\]
\[
=e^{\prime }P_{t}^{-+}P_{t-1}^{-+}...P_{2}^{-+}\left( P_{1}^{K}\right)
^{-+}, 
\]
\noindent \textit{where }$e=e\left( r\right) .$ \textit{Further, since} 
\[
p_{n}^{\prime }=P_{n}P_{n-1}...P_{1}p_{0}^{\prime } 
\]
\noindent \textit{and} (\textit{see Remark} 2.4(c2)) 
\[
\left( P_{n}P_{n-1}...P_{1}\right) ^{K}=\left( P_{t}P_{t-1}...P_{1}\right)
^{K},\text{ }\forall n\geq t,\hspace{0.2cm}\text{\textit{if} }\left|
K^{c}\right| >1 
\]
\noindent \textit{and} 
\[
P_{n}P_{n-1}...P_{1}=P_{t}P_{t-1}...P_{1}\hspace{0.2cm}\text{\textit{and...,}
}\forall n\geq t,\hspace{0.2cm}\text{\textit{if} }\left| K^{c}\right| \leq
1, 
\]

\noindent \textit{we have} 
\[
p_{nl}=\left\{ 
\begin{array}{l}
\sum\limits_{j\in \left\langle r\right\rangle }\left(
P_{t}^{-+}P_{t-1}^{-+}...P_{1}^{-+}\right) _{\left\langle r\right\rangle
\left\{ j\right\} }p_{0j}\smallskip \\ 
\hspace{2cm}\text{\textit{if} }\left| K^{c}\right| \leq 1,\medskip \\ 
\sum\limits_{j\in K}\left( P_{t}^{-+}P_{t-1}^{-+}...P_{1}^{-+}\right)
_{\left\langle r\right\rangle \left\{ j\right\} }p_{0j}+\sum\limits_{u\in
K^{c}}\left( P_{n}P_{n-1}...P_{1}\right) _{lu}p_{0u}\smallskip \\ 
\hspace{2cm}\text{\textit{if} }\left| K^{c}\right| >1,%
\end{array}
\right. 
\]
\noindent $\forall n\geq t,$ $\forall l\in \left\langle r\right\rangle $ (%
\textit{when} $\left| K^{c}\right| >1,$ \textit{we can replace }$\left(
P_{t}^{-+}P_{t-1}^{-+}...P_{1}^{-+}\right) _{\left\langle r\right\rangle
\left\{ j\right\} }$\textit{\ with} $\left(
P_{t}^{-+}P_{t-1}^{-+}...P_{2}^{-+}\left( P_{1}^{K}\right) ^{-+}\right)
_{\left\langle r\right\rangle \left\{ j\right\} };$ $%
P_{t}^{-+}P_{t-1}^{-+}...P_{1}^{-+}$ \textit{has just one row, and the label
of this row is }$\left\langle r\right\rangle $) --- \textit{when} $\emptyset
\neq K\subseteq \left\langle r\right\rangle $ ($K=\left\langle
r\right\rangle $ \textit{or} $K\neq \left\langle r\right\rangle $), \textit{%
the linear combination of }$p_{0j},$\textit{\ }$j\in K,$\textit{\ namely, }$%
\sum\limits_{j\in K}\left( P_{t}^{-+}P_{t-1}^{-+}...P_{1}^{-+}\right)
_{\left\langle r\right\rangle \left\{ j\right\} }p_{0j},$ \textit{does not
depend on }$n$\textit{\ and }$l$\textit{\ when }$n\geq t$\textit{, and,
therefore, a partial or total consensus is reached at time }$t$ (\textit{the
linear combination of }$p_{0j},$\textit{\ }$j\in K,$\textit{\ at time/step }$%
t$\textit{\ is kept at any time }$n>t$\textit{\ --- the utilization of
phrases }``\textit{partial consensus}''\textit{\ and }``\textit{total
consensus}''\textit{\ are thus justified}); \textit{when} $\left|
K^{c}\right| \leq 1,$ \textit{a }(\textit{total})\textit{\ consensus is
reached at time }$t$\textit{\ and this is} 
\[
\sum\limits_{j\in \left\langle r\right\rangle }\left(
P_{t}^{-+}P_{t-1}^{-+}...P_{1}^{-+}\right) _{\left\langle r\right\rangle
\left\{ j\right\} }p_{0j}; 
\]
\noindent \textit{when }$\left| K^{c}\right| >1,$%
\[
\sum\limits_{j\in K}\left( P_{t}^{-+}P_{t-1}^{-+}...P_{1}^{-+}\right)
_{\left\langle r\right\rangle \left\{ j\right\} }p_{0j} 
\]
\noindent \textit{is a partial consensus reached at time }$t$\textit{\ for
the initial opinions, }$p_{0j},$ $j\in K$, \textit{of individuals} \textit{%
from} $K$ \textit{and} \textit{the linear combination of }$p_{0u},$\textit{\ 
}$u\in K^{c},$\textit{\ namely, }$\sum\limits_{u\in K^{c}}\left(
P_{n}P_{n-1}...P_{1}\right) _{lu}p_{0u},$ \textit{depends or not on} $n$ 
\textit{and }$l$. (See also Remark 1.11.)

\smallskip

\textbf{Proof.} Theorems 1.10 and 2.2 and Remarks 1.2(b) and 2.4, (b) or
(c). (Using Remark 1.2(b), we have 
\[
\left( P_{n}P_{n-1}...P_{1}\right) ^{K}=P_{n}P_{n-1}...P_{t+1}\left(
P_{t}P_{t-1}...P_{1}\right) ^{K}=\left( P_{t}P_{t-1}...P_{1}\right) ^{K},%
\text{ }\forall n\geq t 
\]
\noindent --- the latter equation follows from the fact that $\left(
P_{t}P_{t-1}...P_{1}\right) ^{K}$ is a stable (stochastic) matrix ($%
P_{n}P_{n-1}...P_{t+1}$ vanishes when $n\!=\!t$).) 
\endproof%

\smallskip

\textbf{Remark 2.6.} (a) $P\in G_{\left( \left\{ i\right\} \right) _{i\in
\left\langle n\right\rangle },\left( \left\{ i\right\} \right) _{i\in
\left\langle n\right\rangle }},$ $\forall P\in S_{n},$ where $n\geq 1$ (by
Remark 1.7(a)).

\smallskip

(b) A case, an interesting one, of Theorem 2.5 is when $\Delta _{1}=\Delta
_{2}=...=\Delta _{u}=\left( \left\{ i\right\} \right) _{i\in \left\langle
r\right\rangle }$ and $\Delta _{u+1}\neq \left( \left\{ i\right\} \right)
_{i\in \left\langle r\right\rangle }$ for some $u,$ $2\leq u\leq t$ (due to
(a), ``$P_{1}\in G_{\Delta _{2},\Delta _{1}},$ $P_{2}\in G_{\Delta
_{3},\Delta _{2}},...,$ $P_{u-1}\in G_{\Delta _{u},\Delta _{u-1}}$'' holds
without problems). In this case, a (total) consensus is reached.

\smallskip

(c) On Theorem 2.5, our interest is to have a set $K$ as large as possible
and a natural number $t\geq 1$ as small as possible --- interesting
optimization problems.

\smallskip

Below we give two examples to illustrate Theorem 2.5; Remark 2.6(b) is
illustrated in the next example.

\smallskip

\textbf{Example 2.7.} Consider the DeGroot model on $\left( \left\langle
r\right\rangle ,p_{0},\left( P_{n}\right) _{n\geq 1}\right) $, where $r=4,$ $%
p_{0}=\left( 2,4,1,1\right) ,$ and the first four matrices are 
\[
P_{1}=\left( 
\begin{array}{cccc}
\frac{1}{4}\smallskip & 0 & \frac{1}{4} & \frac{2}{4} \\ 
0 & \frac{1}{4}\smallskip & 0 & \frac{3}{4} \\ 
\frac{1}{8} & \frac{1}{8} & \frac{2}{4}\smallskip & \frac{1}{4} \\ 
\frac{3}{16} & \frac{1}{16} & \frac{2}{4} & \frac{1}{4}%
\end{array}
\right) ,\text{ }P_{2}=\left( 
\begin{array}{cccc}
\frac{1}{4}\smallskip & 0 & 0 & \frac{3}{4} \\ 
\frac{1}{4} & 0\smallskip & 0 & \frac{3}{4} \\ 
\frac{3}{4} & 0 & \frac{1}{4}\smallskip & 0 \\ 
\frac{3}{4} & 0 & \frac{1}{4} & 0%
\end{array}
\right) , 
\]
\[
P_{3}=\left( 
\begin{array}{cccc}
\frac{1}{4}\smallskip & 0 & \frac{2}{4} & \frac{1}{4} \\ 
0 & \frac{1}{4}\smallskip & \frac{1}{4} & \frac{2}{4} \\ 
\frac{1}{4} & 0 & 0\smallskip & \frac{3}{4} \\ 
\frac{1}{8} & \frac{1}{8} & \frac{3}{4} & 0%
\end{array}
\right) ,\text{ and }P_{4}=\left( 
\begin{array}{cccc}
\frac{1}{8}\smallskip & \frac{2}{8} & \frac{3}{8} & \frac{2}{8} \\ 
\frac{2}{8} & \frac{3}{8}\smallskip & \frac{2}{8} & \frac{1}{8} \\ 
0 & 0 & 0\smallskip & 1 \\ 
0 & 1 & 0 & 0%
\end{array}
\right) . 
\]
\noindent We have $P_{1}\in G_{\Delta _{2},\Delta _{1}},$ $P_{2}\in
G_{\Delta _{3},\Delta _{2}},$ and $P_{3}\in G_{\Delta _{4},\Delta _{3}},$
where $\Delta _{1}=\Delta _{2}=\left( \left\{ j\right\} \right) _{j\in
\left\langle 4\right\rangle }$ (see now Remark 2.6(b)), $\Delta _{3}=\left(
\left\langle 2\right\rangle ,\left\{ 3,4\right\} \right) ,$ and $\Delta
_{4}=\left( \left\langle 4\right\rangle \right) $ ($P_{1}\in G_{\Delta
_{4},\Delta _{3}},$ but this thing does not count here). By Theorem 2.5 (or
by Theorem 1.10), $P_{3}P_{2}P_{1}$ is a stable matrix ($P_{3}P_{2}$ is also
a stable matrix). Further, we have 
\[
P_{3}^{-+}P_{2}^{-+}P_{1}^{-+}=\left( 
\begin{array}{cc}
\frac{1}{4} & \frac{3}{4}%
\end{array}
\right) \left( 
\begin{array}{cccc}
\frac{1}{4}\smallskip & 0 & 0 & \frac{3}{4} \\ 
\frac{3}{4} & 0 & \frac{1}{4} & 0%
\end{array}
\right) \left( 
\begin{array}{cccc}
\frac{1}{4}\smallskip & 0 & \frac{1}{4} & \frac{2}{4} \\ 
0 & \frac{1}{4}\smallskip & 0 & \frac{3}{4} \\ 
\frac{1}{8} & \frac{1}{8} & \frac{2}{4}\smallskip & \frac{1}{4} \\ 
\frac{3}{16} & \frac{1}{16} & \frac{2}{4} & \frac{1}{4}%
\end{array}
\right) = 
\]
\[
=\left( 
\begin{array}{cccc}
\frac{55}{256} & \frac{9}{256} & \frac{88}{256} & \frac{104}{256}%
\end{array}
\right) , 
\]
\noindent so, the consensus is 
\[
P_{3}^{-+}P_{2}^{-+}P_{1}^{-+}p_{0}^{\prime }=\frac{110+36+88+104}{256}=%
\frac{338}{256} 
\]
\noindent (it is reached at time $3$). The sequence of matrices $\left(
P_{n}\right) _{n\geq 4}$ does not count --- this can be replaced with any
other infinite sequence of stochastic $4\times 4$ matrices (see Remark
2.4(c1)).

\smallskip

\textbf{Example 2.8.} Consider the DeGroot model on $\left( \left\langle
r\right\rangle ,p_{0},\left( P_{n}\right) _{n\geq 1}\right) ,$ where $r=4,$ $%
p_{0}=\left( 10,20,2,4\right) ,$ and the first two matrices are 
\[
P_{1}=\left( 
\begin{array}{cccc}
\frac{2}{10}\smallskip & \frac{1}{10} & \frac{7}{10} & 0 \\ 
\frac{2}{10} & \frac{1}{10}\smallskip & \frac{1}{10} & \frac{6}{10} \\ 
0 & \frac{3}{10} & 0\smallskip & \frac{7}{10} \\ 
0 & \frac{3}{10} & \frac{2}{10} & \frac{5}{10}%
\end{array}
\right) \text{ and }P_{2}=\left( 
\begin{array}{cccc}
\frac{4}{10}\smallskip & 0 & 0 & \frac{6}{10} \\ 
\frac{1}{10} & \frac{3}{10}\smallskip & \frac{4}{10} & \frac{2}{10} \\ 
\frac{1}{10} & \frac{3}{10} & \frac{1}{10}\smallskip & \frac{5}{10} \\ 
\frac{2}{10} & \frac{2}{10} & \frac{1}{10} & \frac{5}{10}%
\end{array}
\right) . 
\]
\noindent We have $P_{1}\in G_{\Delta _{2},\Delta _{1}}$ and $P_{2}\in
G_{\Delta _{3},\Delta _{2}},$ where $\Delta _{1}=\left( \left\{ 1\right\}
,\left\{ 2\right\} ,\left\{ 3,4\right\} \right) ,$ $\Delta _{2}=\left(
\left\{ 1,2\right\} ,\left\{ 3,4\right\} \right) $, and $\Delta _{3}=\left(
\left\langle 4\right\rangle \right) .$ So, by Theorem 2.5 (or by Theorem
1.10), 
\[
\left( P_{2}P_{1}\right) ^{\left\{ 1,2\right\} } 
\]
\noindent is a stable matrix. Further, we have 
\[
P_{2}^{-+}P_{1}^{-+}=\left( 
\begin{array}{cc}
\frac{4}{10} & \frac{6}{10}%
\end{array}
\right) \left( 
\begin{array}{ccc}
\frac{2}{10}\smallskip & \frac{1}{10} & \frac{7}{10} \\ 
0 & \frac{3}{10} & \frac{7}{10}%
\end{array}
\right) =\left( 
\begin{array}{ccc}
\frac{8}{100} & \frac{22}{100} & \frac{70}{100}%
\end{array}
\right) . 
\]
\noindent So (see Theorem 2.5) 
\[
p_{nl}=\frac{8}{100}\cdot 10+\frac{22}{100}\cdot 20+\left(
P_{n}P_{n-1`}...P_{1}\right) _{l3}\cdot 2+\left(
P_{n}P_{n-1`}...P_{1}\right) _{l4}\cdot 4= 
\]
\[
=5.2+2\left( P_{n}P_{n-1`}...P_{1}\right) _{l3}+4\left(
P_{n}P_{n-1`}...P_{1}\right) _{l4},\text{ }\forall n\geq 2,\text{ }\forall
l\in \left\langle 4\right\rangle 
\]
\noindent --- the linear combination, 
\[
\frac{8}{100}\cdot 10+\frac{22}{100}\cdot 20, 
\]
\noindent for the initial opinions $10$ and $20$ (of individuals $1$ and $2$%
) holds at any time $n\geq 2;$ the consensus, this is a partial consensus,
for the initial opinions $10$ and $20$ is $5.2,$ and is reached at time $2.$

\smallskip

Below we consider a possible case of the DeGroot model on $\left(
\left\langle r\right\rangle ,p_{0},\left( P_{n}\right) _{n\geq 1}\right) $
--- a submodel of the DeGroot model which works from small subsets
(subgroups of individuals) of $\left\langle r\right\rangle $ to larger and
larger subsets of $\left\langle r\right\rangle $ until the whole set $%
\left\langle r\right\rangle .$

\smallskip

Let $\Delta _{1},\Delta _{2},...,\Delta _{t+1}\in $Par$\left( \left\langle
r\right\rangle \right) ,$ $r\geq 2,$ $t\geq 1.$ Suppose that $\Delta
_{1}\preceq \Delta _{2}\preceq ...\preceq \Delta _{t+1},$ $\Delta
_{1}=\left( \left\{ j\right\} \right) _{j\in \left\langle r\right\rangle },$
and $\Delta _{t+1}=\left( \left\langle r\right\rangle \right) .$ Consider
the DeGroot model on $\left( \left\langle r\right\rangle ,p_{0},\left(
P_{n}\right) _{n\geq 1}\right) ,$ where the matrices $%
P_{1},P_{2},...,P_{t-1} $ (not $P_{1},P_{2},...,P_{t}$) when $t\geq 2$ have
the property 
\[
\left( P_{l}\right) _{K}^{L}=0,\text{ }\forall l\in \left\langle
t-1\right\rangle \text{, }\forall K,L\in \Delta _{l+1},\text{ }K\neq L 
\]
\noindent (this assumption implies that $P_{l}$ is a block diagonal matrix,
eventually by permutation of rows and columns, and $\Delta _{l+1}$-stable
matrix on $\Delta _{l+1},$ $\forall l\in \left\langle t-1\right\rangle $).

\smallskip

\textbf{Remark 2.9.} Theorem 2.5 holds, in particular, for the above DeGroot
submodel when $P_{1}\in G_{\Delta _{2},\Delta _{1}},$ $P_{2}\in G_{\Delta
_{3},\Delta _{2}},$ $...,$ $P_{t}\in G_{\Delta _{t+1},\Delta _{t}}.$

\smallskip

To illustrate the DeGroot submodel, we give two examples.

\smallskip

\textbf{Example 2.10. }Consider the DeGroot model on $\left( \left\langle
r\right\rangle ,p_{0},\left( P_{n}\right) _{n\geq 1}\right) ,$ where 
\[
r=6,\text{ }p_{0}=\left( \frac{1}{10},\frac{2}{10},\frac{3}{10},\frac{1}{10},%
\frac{1}{10},\frac{2}{10}\right) , 
\]
\noindent and the first three matrices are 
\[
P_{1}=\left( 
\begin{array}{cccccc}
\frac{1}{2}\smallskip & \frac{1}{2} &  &  &  &  \\ 
\frac{2}{3} & \frac{1}{3}\smallskip &  &  &  &  \\ 
&  & \frac{1}{4}\smallskip & \frac{3}{4} &  &  \\ 
&  & \frac{2}{4} & \frac{2}{4}\smallskip &  &  \\ 
&  &  &  & 0\smallskip & 1 \\ 
&  &  &  & 1 & 0%
\end{array}
\right) ,\text{ }P_{2}=\left( 
\begin{array}{cccccc}
1\smallskip & 0 & 0 & 0 &  &  \\ 
\frac{1}{2} & 0\smallskip & \frac{1}{2} & 0 &  &  \\ 
\frac{1}{4} & \frac{1}{4} & \frac{1}{4}\smallskip & \frac{1}{4} &  &  \\ 
0 & 1 & 0 & 0\smallskip &  &  \\ 
&  &  &  & 1\smallskip & 0 \\ 
&  &  &  & \frac{1}{10} & \frac{9}{10}%
\end{array}
\right) , 
\]
\noindent and 
\[
P_{3}=\left( 
\begin{array}{cccccc}
1\smallskip & 0 & 0 & 0 & 0 & 0 \\ 
\frac{1}{6} & \frac{1}{6}\smallskip & \frac{1}{6} & \frac{1}{6} & \frac{1}{6}
& \frac{1}{6} \\ 
0 & 1 & 0\smallskip & 0 & 0 & 0 \\ 
0 & 1 & 0 & 0\smallskip & 0 & 0 \\ 
0 & 0 & 1 & 0 & 0\smallskip & 0 \\ 
\frac{1}{6} & \frac{1}{6} & \frac{1}{6} & \frac{1}{6} & \frac{1}{6} & \frac{1%
}{6}%
\end{array}
\right) . 
\]
\noindent Taking $\Delta _{1}=\left( \left\{ j\right\} \right) _{j\in
\left\langle 6\right\rangle },$ $\Delta _{2}=\left( \left\langle
2\right\rangle ,\left\{ 3,4\right\} ,\left\{ 5,6\right\} \right) ,$ $\Delta
_{3}=\left( \left\langle 4\right\rangle ,\left\{ 5,6\right\} \right) ,$ and $%
\Delta _{4}=\left( \left\langle 6\right\rangle \right) ,$ we have $\Delta
_{1}\preceq \Delta _{2}\preceq \Delta _{3}\preceq \Delta _{4}$ and 
\[
\left( P_{1}\right) _{\left\langle 2\right\rangle }^{\left\{ 3,4\right\}
}\equiv \left( P_{1}\right) _{\left\langle 2\right\rangle }^{\left\{
5,6\right\} }=0,\text{ }\left( P_{1}\right) _{\left\{ 3,4\right\}
}^{\left\langle 2\right\rangle }\equiv \left( P_{1}\right) _{\left\{
3,4\right\} }^{\left\{ 5,6\right\} }=0, 
\]
\[
\left( P_{1}\right) _{\left\{ 5,6\right\} }^{\left\langle 2\right\rangle
}\equiv \left( P_{1}\right) _{\left\{ 5,6\right\} }^{\left\{ 3,4\right\} }=0,%
\text{ }\left( P_{2}\right) _{\left\langle 4\right\rangle }^{\left\{
5,6\right\} }=0,\text{ and }\left( P_{2}\right) _{\left\{ 5,6\right\}
}^{\left\langle 4\right\rangle }=0 
\]
\noindent --- the conditions of the DeGroot submodel are satisfied, and thus
we have an example for this submodel. Learning begin with the small subsets
(subgroups of individuals) $\left\langle 2\right\rangle ,$ $\left\{
3,4\right\} ,$ and $\left\{ 5,6\right\} $ of $\left\langle 6\right\rangle $,
the larger subsets $\left\langle 4\right\rangle $ and $\left\{ 5,6\right\} $
of $\left\langle 6\right\rangle $ are then used, then the whole set, $%
\left\langle 6\right\rangle ,$ is used. Theorem 2.5 cannot be applied for $%
\Delta _{1},$ $\Delta _{2},$ $\Delta _{3},$ and $\Delta _{4}$ because $%
P_{1}\notin G_{\Delta _{2},\Delta _{1}}$ (moreover, $P_{2}\notin G_{\Delta
_{3},\Delta _{2}}$).

\smallskip

In the next example we show that our DeGroot submodel on$\left( \left\langle
r\right\rangle \!,p_{0},\!\left( P_{n}\right) _{n\geq 1}\right) $ has a
subset/subgroup property when $P_{1}\in G_{\Delta _{2},\Delta _{1}},$ $%
P_{2}\in G_{\Delta _{3},\Delta _{2}},$ $...,$ $P_{t}\in G_{\Delta
_{t+1},\Delta _{t}}.$

\smallskip

\textbf{Example 2.11.} (The subset/subgroup property illustrated.) Consider
the DeGroot model on $\left( \left\langle r\right\rangle ,p_{0},\left(
P_{n}\right) _{n\geq 1}\right) ,$ where 
\[
r=6,\text{ }p_{0}=\left( \frac{1}{10},\frac{2}{10},\frac{3}{10},\frac{1}{10},%
\frac{1}{10},\frac{2}{10}\right) , 
\]
\noindent and the first four matrices are 
\[
P_{1}=\left( 
\begin{array}{cccccc}
\frac{4}{10}\smallskip & \frac{6}{10} &  &  &  &  \\ 
\frac{4}{10} & \frac{6}{10}\smallskip &  &  &  &  \\ 
&  & 1\smallskip &  &  &  \\ 
&  &  & 1\smallskip &  &  \\ 
&  &  &  & \frac{1}{10}\smallskip & \frac{9}{10} \\ 
&  &  &  & \frac{1}{10} & \frac{9}{10}%
\end{array}
\right) ,\text{ }P_{2}=\left( 
\begin{array}{cccccc}
\frac{3}{10}\smallskip & \frac{4}{10} & \frac{3}{10} &  &  &  \\ 
\frac{2}{10} & \frac{5}{10}\smallskip & \frac{3}{10} &  &  &  \\ 
\frac{1}{10} & \frac{6}{10} & \frac{3}{10}\smallskip &  &  &  \\ 
&  &  & 1\smallskip &  &  \\ 
&  &  &  & \frac{2}{10}\smallskip & \frac{8}{10} \\ 
&  &  &  & \frac{4}{10} & \frac{6}{10}%
\end{array}
\right) , 
\]
\[
P_{3}=\left( 
\begin{array}{cccccc}
\frac{3}{10}\smallskip & \frac{2}{10} & \frac{5}{10} &  &  &  \\ 
\frac{2}{10} & \frac{4}{10}\smallskip & \frac{4}{10} &  &  &  \\ 
\frac{1}{10} & \frac{6}{10} & \frac{3}{10}\smallskip &  &  &  \\ 
&  &  & \frac{1}{10}\smallskip & \frac{1}{10} & \frac{8}{10} \\ 
&  &  & \frac{1}{10} & \frac{2}{10}\smallskip & \frac{7}{10} \\ 
&  &  & \frac{1}{10} & \frac{1}{10} & \frac{8}{10}%
\end{array}
\right) ,\text{ }P_{4}=\left( 
\begin{array}{cccccc}
\frac{2}{10}\smallskip & \frac{2}{10} & \frac{2}{10} & \frac{1}{10} & \frac{1%
}{10} & \frac{2}{10} \\ 
\frac{1}{10} & \frac{3}{10}\smallskip & \frac{2}{10} & \frac{1}{10} & \frac{2%
}{10} & \frac{1}{10} \\ 
\frac{1}{10} & \frac{2}{10} & \frac{3}{10}\smallskip & \frac{2}{10} & \frac{1%
}{10} & \frac{1}{10} \\ 
\frac{3}{10} & \frac{2}{10} & \frac{1}{10} & \frac{1}{10}\smallskip & \frac{1%
}{10} & \frac{2}{10} \\ 
\frac{1}{10} & \frac{4}{10} & \frac{1}{10} & \frac{2}{10} & \frac{1}{10}%
\smallskip & \frac{1}{10} \\ 
\frac{1}{10} & \frac{4}{10} & \frac{1}{10} & \frac{1}{10} & \frac{1}{10} & 
\frac{2}{10}%
\end{array}
\right) . 
\]
\noindent It is easy to see that this model satisfies the conditions of
DeGroot submodel for $\Delta _{1}=\left( \left\{ j\right\} \right) _{j\in
\left\langle 6\right\rangle },$ $\Delta _{2}=\left( \left\langle
2\right\rangle ,\left\{ 3\right\} ,\left\{ 4\right\} ,\left\{ 5,6\right\}
\right) ,$ $\Delta _{3}=\left( \left\langle 3\right\rangle ,\left\{
4\right\} ,\left\{ 5,6\right\} \right) ,$ $\Delta _{4}=\left( \left\langle
3\right\rangle ,\left\{ 4,5,6\right\} \right) ,$ and $\Delta _{5}=\left(
\left\langle 6\right\rangle \right) .$ We have $P_{1}\in G_{\Delta
_{2},\Delta _{1}},$ $P_{2}\in G_{\Delta _{3},\Delta _{2}},$ $P_{3}\in
G_{\Delta _{4},\Delta _{3}},$ and $P_{4}\in G_{\Delta _{5},\Delta _{4}}.$ By
Theorem 2.5, $P_{4}P_{3}P_{2}P_{1}$ is a stable matrix, so, a consensus is
reached at time $4$... --- the continuation, using Theorem 2.5, is left to
the reader (see also Example 2.7). Further, we comment on the matrices $%
P_{1},$ $P_{2},$ $P_{3},$ and $P_{4},$ and use Theorem 1.16. An
interpretation of $P_{1}$ (at time $1$ --- see the DeGroot model again): the
individuals $1$ and $2$ negotiate --- negotiation, a good idea ---, and
agree the same probability, $\frac{4}{10},$ for the initial
opinion/estimation of individual $1$, \textit{i.e.}, for $p_{01},$ and,
separately, the same probability, $\frac{6}{10},$ for the initial opinion of
individual $2$, \textit{i.e.}, for $p_{02};$ the individuals $5$ and $6$
proceed similarly; the individuals $3$ and $4$ do not negotiate... An
interpretation of $P_{2}$ (at time $2$): the individuals $1$, $2$, and $3$
negotiate, and agree the same probability, $\frac{7}{10}$ ($\frac{7}{10}=%
\frac{3}{10}+\frac{4}{10}=\frac{2}{10}+\frac{5}{10}=\frac{1}{10}+\frac{6}{10}
$), for the opinions of individuals $1$ and $2$ at time $1$ considered
together, \textit{i.e.}, for $p_{11}$ and $p_{12}$ considered together, and,
separately, the same probability, $\frac{3}{10},$ for the opinion of
individual $3$ at time $1$, \textit{i.e.}, for $p_{13};$ the individuals $5$
and $6$ renegotiate, and agree the same probability, 1 ($1=\frac{2}{10}+%
\frac{8}{10}=\frac{4}{10}+\frac{6}{10}$), for their opinions at time $1$
considered together, \textit{i.e.}, for $p_{15}$ and $p_{16}$ considered
together; the individual $4$ does not negotiate... An interpretation of $%
P_{3}$ (at time $3$): the individuals $1$, $2$, and $3$ renegotiate, and
agree the same probability, $1$ ($1=\frac{3}{10}+\frac{2}{10}+\frac{5}{10}=%
\frac{2}{10}+\frac{4}{10}+\frac{4}{10}=\frac{1}{10}+\frac{6}{10}+\frac{3}{10}
$), for their opinions at time $2$ considered together, \textit{i.e.}, for $%
p_{21},$ $p_{22},$ and $p_{23}$ considered together; the individuals $4,$ $%
5, $ and $6$ negotiate, and agree the same probability, $\frac{1}{10}$, for
the opinion of individual $4$ at time $2$, \textit{i.e.}, for $p_{24},$ and,
separately, the same probability, $\frac{9}{10}$ ($\frac{9}{10}=\frac{1}{10}+%
\frac{8}{10}=\frac{2}{10}+\frac{7}{10}=\frac{1}{10}+\frac{8}{10}$), for the
opinions of individuals $5$ and $6$ at time $2$ considered together, \textit{%
i.e.}, for $p_{25}$ and $p_{26}$ considered together. An interpretation of $%
P_{4}$ (at time $4$): all the individuals negotiate, and agree the same
probability, $\frac{6}{10}$ ($\frac{6}{10}=\frac{2}{10}+\frac{2}{10}+\frac{2%
}{10}=\frac{1}{10}+\frac{3}{10}+\frac{2}{10}=...$), for the opinions of
individuals $1$, $2$, and $3$ at time $3$ considered together, \textit{i.e.}%
, for $p_{31},$ $p_{32},$ and $p_{33}$ considered together, and, separately,
the same probability, $\frac{4}{10}$ ($\frac{4}{10}=\frac{1}{10}+\frac{1}{10}%
+\frac{2}{10}=\frac{1}{10}+\frac{2}{10}+\frac{1}{10}=...$), for the opinions
of individuals $4$, $5$, and $6$ at time $3$ considered together, \textit{%
i.e.}, for $p_{34},$ $p_{35},$ and $p_{36}$ considered together. The
negotiations and renegotiations took place from small subsets (of $%
\left\langle 6\right\rangle $)/subgroups to larger and larger
subsets/subgroups until to the whole set/group. By Theorem 1.16, $P_{2}$ can
be replaced with any matrix which is similar to it, such as, 
\[
\left( 
\begin{array}{cccccc}
\frac{3}{10}\smallskip & \frac{4}{10} & \frac{3}{10} &  &  &  \\ 
\frac{3}{10} & \frac{4}{10}\smallskip & \frac{3}{10} &  &  &  \\ 
\frac{3}{10} & \frac{4}{10} & \frac{3}{10}\smallskip &  &  &  \\ 
&  &  & 1\smallskip &  &  \\ 
&  &  &  & 1\smallskip & 0 \\ 
&  &  &  & 0 & 1%
\end{array}
\right) :=P_{2}^{\prime }, 
\]
\noindent and 
\[
\left( 
\begin{array}{cccccc}
\frac{7}{10}\smallskip & 0 & \frac{3}{10} &  &  &  \\ 
0 & \frac{7}{10}\smallskip & \frac{3}{10} &  &  &  \\ 
\frac{7}{10} & 0 & \frac{3}{10}\smallskip &  &  &  \\ 
&  &  & 1\smallskip &  &  \\ 
&  &  &  & 0\smallskip & 1 \\ 
&  &  &  & 1 & 0%
\end{array}
\right) :=P_{2}^{\prime \prime }. 
\]
\noindent The matrices $\left( P_{2}\right) _{\left\langle 3\right\rangle
}^{\left\langle 2\right\rangle },$%
\[
\left( P_{2}\right) _{\left\langle 3\right\rangle }^{\left\langle
2\right\rangle }=\left( 
\begin{array}{cc}
\frac{3}{10}\smallskip & \frac{4}{10} \\ 
\frac{2}{10} & \frac{5}{10}\smallskip \\ 
\frac{1}{10} & \frac{6}{10}%
\end{array}
\right) , 
\]
\noindent $\left( P_{2}^{\prime }\right) _{\left\langle 3\right\rangle
}^{\left\langle 2\right\rangle },$%
\[
\left( P_{2}^{\prime }\right) _{\left\langle 3\right\rangle }^{\left\langle
2\right\rangle }=\left( 
\begin{array}{cc}
\frac{3}{10}\smallskip & \frac{4}{10} \\ 
\frac{3}{10} & \frac{4}{10}\smallskip \\ 
\frac{3}{10} & \frac{4}{10}%
\end{array}
\right) , 
\]
\noindent and $\left( P_{2}^{\prime \prime }\right) _{\left\langle
3\right\rangle }^{\left\langle 2\right\rangle },$%
\[
\left( P_{2}^{\prime \prime }\right) _{\left\langle 3\right\rangle
}^{\left\langle 2\right\rangle }=\left( 
\begin{array}{cc}
\frac{7}{10}\smallskip & 0 \\ 
0 & \frac{7}{10}\smallskip \\ 
\frac{7}{10} & 0%
\end{array}
\right) 
\]
\noindent (each of them has $6$ entries, 
\[
\left( 
\begin{array}{cc}
\star & \star \\ 
\star & \star \\ 
\star & \star%
\end{array}
\right) , 
\]
\noindent ``$\star $'' stands for a positive or zero entry), have the
property, a subset/subgroup property: the row sums are equal to $\frac{7}{10}
$ --- only the row sums count, not the positions of positive entries, but at
least one entry must be positive in each row (because $\frac{7}{10}>0$),
only $\frac{7}{10}$ counts, $\frac{7}{10}$ is the probability... (see
above), and this probability is for the subset $\left\langle 2\right\rangle $
of $\left\langle 6\right\rangle $. Among $P_{2},$ $P_{2}^{\prime },$ and $%
P_{2}^{\prime \prime },$ $P_{2}^{\prime \prime }$ is the best because it has
the smallest number of positive entries --- these are the positive
probabilities for the subsets involved here: $\frac{7}{10}$ for (the subset) 
$\left\langle 2\right\rangle $ ($\frac{7}{10}=$ the probability assigned by
the individual $i$ for the opinions of individuals 1 and 2 at time 1
considered together, $\forall i\in \left\langle 3\right\rangle $), $\frac{3}{%
10}$ for $\left\{ 3\right\} $, $1$ for $\left\{ 4\right\} $, and $1$ for $%
\left\{ 5,6\right\} $. Similarly, $P_{3}$ and $P_{4}$ can be replaced with
any matrices, $Q_{3}$ and $Q_{4}$, which are similar to $P_{3}$ and $P_{4},$
respectively. By Theorem 1.16 we have 
\[
P_{4}P_{3}P_{2}P_{1}=P_{4}P_{3}P_{2}^{\prime }P_{1}=P_{4}P_{3}P_{2}^{\prime
\prime }P_{1}=P_{4}Q_{3}P_{2}P_{1}=... 
\]
\noindent (we can work with $P_{4}P_{3}P_{2}P_{1}$, or with $%
P_{4}P_{3}P_{2}^{\prime }P_{1},$ or with...). Note that $P_{1},$ $P_{2},$
and $P_{3}$ are reducible while $P_{4}$ is irreducible (being positive), but
the last matrix can be replaced with any reducible matrix which is similar
to it, such as, 
\[
\left( 
\begin{array}{cccccc}
\frac{6}{10}\smallskip & 0 & 0 & \frac{4}{10} & 0 & 0 \\ 
\frac{6}{10} & 0\smallskip & 0 & \frac{4}{10} & 0 & 0 \\ 
\frac{6}{10} & 0 & 0\smallskip & \frac{4}{10} & 0 & 0 \\ 
\frac{6}{10} & 0 & 0 & \frac{4}{10}\smallskip & 0 & 0 \\ 
\frac{6}{10} & 0 & 0 & \frac{4}{10} & 0\smallskip & 0 \\ 
\frac{6}{10} & 0 & 0 & \frac{4}{10} & 0 & 0%
\end{array}
\right) :=U_{4} 
\]
\noindent --- $P_{1},$ $P_{2},$ $P_{3},$ and $U_{4}$ are reducible, but,
interestingly, $U_{4}P_{3}P_{2}P_{1}$ is irreducible (because it is positive
--- it is also stable).

\smallskip

For other examples for the DeGroot submodel, see, in the next section,
Remark 3.7.

\bigskip

\bigskip

\bigskip

\begin{center}
\textbf{3. FINITE-TIME\ CONSENSUS\ FOR\ DEGROOT MODEL\ ON\ DISTRIBUTED\
SYSTEMS}
\end{center}

\bigskip

In this section, we consider the finite-time consensus problem for the
DeGroot model on distributed systems (for these systems, see, \textit{e.g.},
[6], see, \textit{e.g.}, also [4]) --- this section can also be useful for
the study of other systems, such as, the multi-agent systems (for the
multi-agent systems, see, \textit{e.g.}, [14], [23], and [29]). For the
DeGroot model on distributed systems, using the $G$ method too, we have a
result for reaching a partial or total (distributed) consensus in a finite
time similar to that for the DeGroot model for reaching a partial or total
consensus in a finite time (see Remark 3.2(b)). We show that for any
connected graph having $2^{m}$ vertices, $m\geq 1,$ and a spanning subgraph
isomorphic to the $m$-cube graph, distributed averaging is performed in $m$
steps --- this result can be extended --- research work --- for any graph
with $n_{1}n_{2}...n_{t}$ vertices under certain conditions, where $t,$ $%
n_{1},n_{2},...,n_{t}\geq 2,$ and, in this case, distributed averaging is
performed in $t$ steps.

\smallskip

Let $\mathcal{G}=\left( \mathcal{V},\mathcal{E}\right) $ be a (nondirected
simple finite) connected graph, where $\mathcal{V}$ is the vertex/node set
and $\mathcal{E}$ is the edge set. (A \textit{simple graph} is a graph
without multiple edges and loops.) Suppose that $\mathcal{V}=\left\langle
r\right\rangle ,$ and that $r\geq 2.$ The edge set $\mathcal{E}$ is a set of
2-element subsets of $\left\langle r\right\rangle ,$ so, 
\[
\mathcal{E}\subseteq \left\{ \left\{ i,j\right\} \left| \text{ }i,j\in
\left\langle r\right\rangle ,\text{ }i\neq j\right. \right\} \text{ (}%
\mathcal{E}\neq \emptyset \text{)} 
\]
\noindent --- \textit{e.g.}, $\mathcal{E}=\left\{ \left\{ 1,2\right\}
,\left\{ 1,3\right\} ,...,\left\{ 1,r\right\} \right\} $ ($r\geq 2$). This
graph, on the distributed systems or not, can be a mathematical model for a
network of computers, processors, or sensors, for a team, when the graphs
are involved, of robots, aircraft, or spacecraft... Set 
\[
\mathcal{N}_{i}=\left\{ j\left| \text{ }j\in \left\langle r\right\rangle 
\text{ and }\left\{ i,j\right\} \in \mathcal{E}\right. \right\} ,\text{ }%
\forall i\in \left\langle r\right\rangle . 
\]
\noindent $\mathcal{N}_{i}$ is called the \textit{set of neighbors of }$i,$ $%
\forall i\in \left\langle r\right\rangle .$ $i\notin \mathcal{N}_{i},$ $%
\forall i\in \left\langle r\right\rangle .$ Suppose that each vertex/node $%
i\in \left\langle r\right\rangle $ holds an initial value, say, $q_{0i},$ $%
q_{0i}\in \mathbb{R},$ and we call the vector 
\[
q_{0}=\left( q_{0i}\right) _{i\in \left\langle r\right\rangle }=\left(
q_{01},q_{02},...,q_{0r}\right) 
\]
\noindent the \textit{initial value vector}. We use homogeneous and
nonhomogeneous linear iterations: at each time $n\geq 1$ ($n\in \mathbb{N}%
^{*}$), each vertex $i$ updates its value as (see, \textit{e.g.}, [27] and
[30] --- only the homogeneous linear iterations are considered in these
articles) 
\[
q_{ni}=\sum\limits_{j\in \left\{ i\right\} \cup \mathcal{N}_{i}}\left(
W_{n}\right) _{ij}q_{n-1j}, 
\]
\noindent where $\left( W_{n}\right) _{ij}\in \mathbb{R},$ $\forall j\in
\left\{ i\right\} \cup \mathcal{N}_{i},$ $\left( W_{n}\right) _{ij}$ is the
weight on $q_{n-1j}$ at vertex $i$ at time $n$, $\forall j\in \left\{
i\right\} \cup \mathcal{N}_{i}$. Set the (row) vector 
\[
q_{n}=\left( q_{ni}\right) _{i\in \left\langle r\right\rangle }=\left(
q_{n1},q_{n2},...,q_{nr}\right) ,\text{ }\forall n\geq 0 
\]
\noindent --- using $q_{0},$ we obtain $q_{1}$ (see above), using $q_{1},$
we obtain $q_{2},$ ... We call the vector $q_{n}$ the \textit{value vector
at time/step }$n,$ $\forall n\geq 0$ --- the value vector at time $0$ is
also (see above) called the initial value vector. Set the real $r\times r$
matrix 
\[
W_{n}=\left( \left( W_{n}\right) _{ij}\right) _{i,j\in \left\langle
r\right\rangle }, 
\]
\noindent $\left( W_{n}\right) _{ij}$ was defined above, $\forall i\in
\left\langle r\right\rangle ,$ $\forall j\in \left\{ i\right\} \cup \mathcal{%
N}_{i},$ and 
\[
\left( W_{n}\right) _{ij}=0\text{ if }j\notin \left\{ i\right\} \cup 
\mathcal{N}_{i},\text{ }\forall i\in \left\langle r\right\rangle 
\]
\noindent (we can have $\left( W_{n}\right) _{ij}=0$ even when $j\in \left\{
i\right\} \cup \mathcal{N}_{i}$, where $i\in \left\langle r\right\rangle $
and $n\geq 1$), $\forall n\geq 1$. We call the matrix $W_{n}$ the \textit{%
weight matrix at time/step} $n,$ $\forall n\geq 1.$ $W_{n},$ $n$ fixed, is
symmetric or not (each (nondirected/undirected) edge can be considered as
two directed edges). Obviously, we have 
\[
q_{n}^{\prime }=W_{n}q_{n-1}^{\prime }=W_{n}W_{n-1}...W_{1}q_{0}^{\prime },%
\text{ }\forall n\geq 1 
\]
\noindent ($q_{n}^{\prime }$ is the transpose of $q_{n}$...).

\smallskip

The above model is a DeGroot-type model --- compare this model with the
DeGroot model from Section 2. It is determined by the graph $\mathcal{G}$,
real vector $q_{0},$ and real matrices $W_{n},$ $n\geq 1.$ Due these facts,
we call it the \textit{DeGroot model on} (\textit{the triple}) $\left( 
\mathcal{G},q_{0},\left( W_{n}\right) _{n\geq 1}\right) ,$ and will use the
generic name ``(the) DeGroot model on distributed systems''. We say that the
DeGroot model on $\left( \mathcal{G},q_{0},\left( W_{n}\right) _{n\geq
1}\right) $ is \textit{homogeneous} if $W_{1}=W_{2}=W_{3}=...$ and \textit{%
nonhomogeneous} if $\exists u,v\geq 1,$ $u\neq v,$ such that $W_{u}\neq
W_{v},$ and will use the generic names ``(the) homogeneous DeGroot model on
distributed systems'' and ``(the) nonhomogeneous DeGroot model on
distributed systems'', respectively.

\smallskip

Note that the DeGroot model on $\left( \left\langle r\right\rangle
,p_{0},\left( P_{n}\right) _{n\geq 1}\right) $ (see Section 2) is trivially
identical with the DeGroot model on $\left( \mathcal{K}_{r},p_{0},\left(
P_{n}\right) _{n\geq 1}\right) ,$ $\forall r\geq 2,$ $\forall p_{0},$ $p_{0}$
is a real vector here, $\forall \left( P_{n}\right) _{n\geq 1},$ $\left(
P_{n}\right) _{n\geq 1}$ is a sequence of stochastic matrices here, where $%
\mathcal{K}_{r}$ is the complete graph with $r$ vertices, $\forall r\geq 2$
--- moreover, it is identical with the DeGroot model on $\left( \mathcal{G}%
_{1},p_{0},\left( P_{n}\right) _{n\geq 1}\right) ,$ $\forall r\geq 2,$ $%
\forall p_{0},$ $p_{0}$ is a real vector here, $\forall \left( P_{n}\right)
_{n\geq 1},$ $\left( P_{n}\right) _{n\geq 1}$ is a sequence of stochastic
matrices here, where $\mathcal{G}_{1}=\left( \mathcal{V}_{1},\mathcal{E}%
_{1}\right) ,$ $\mathcal{V}_{1}=\left\langle r\right\rangle $ and

$\!\!\!\!\!\!\mathcal{E}_{1}\!=\!\left\{ \left\{ i,j\right\} \left| \text{ }%
i,j\in \left\langle r\right\rangle \!,\text{ }\!i\neq j,\!\text{ and }%
\exists n\!\geq \!1\text{ such that }\left( P_{n}\right) _{ij}\!>\!0\text{
or }\left( P_{n}\right) _{ji}\!>\!0\right. \!\right\} \!\!.$

\smallskip

Definition 2.1 (from Section 2) and the first part of the next definition
are similar.

\smallskip

\textbf{Definition 3.1.} (See also [27] and [30].) Consider the DeGroot
model on $\left( \mathcal{G},q_{0},\left( W_{n}\right) _{n\geq 1}\right) .$
We say that the graph $\mathcal{G}$ reaches a (\textit{distributed}) \textit{%
consensus} if $\lim\limits_{n\rightarrow \infty }W_{n}W_{n-1}...W_{1}$
exists and is a stable matrix. Set 
\[
W_{\infty }=\lim\limits_{n\rightarrow \infty }W_{n}W_{n-1}...W_{1} 
\]
\noindent when $\lim\limits_{n\rightarrow \infty }W_{n}W_{n-1}...W_{1}$
exists. ($W_{\infty }$ is a real $r\times r$ matrix.) In particular, if 
\[
W_{\infty }=e^{\prime }\left( \frac{1}{r},\frac{1}{r},...,\frac{1}{r}\right)
, 
\]
\noindent we say that the graph $\mathcal{G}$ \textit{performs distributed
averaging} ($e\!=\!e\left( r\right) \!=\!\left( 1,1,...,1\right) ;$ $%
e^{\prime }$ is its transpose).

\smallskip

Reaching a distributed consensus for the DeGroot model on distributed
systems --- this is the (\textit{distributed})\textit{\ consensus problem}
for/of the DeGroot model on distributed systems, and is similar to the
consensus problem for the DeGroot model. When a distributed consensus is
reached, it is reached in a finite or infinite time --- obviously, the first
case is the best. When distributed averaging is performed, it is performed
in a finite or infinite time --- the first case is also the best.

\smallskip

When the graph $\mathcal{G}$ performs distributed averaging, each vertex
holds the average of initial values, $\frac{1}{r}\sum\limits_{i=1}^{r}q_{0i}$
--- an interpretation (see, \textit{e.g.}, [30, p. 65]): if $q_{0i}$ is the
amount of a resource at the vertex $i,$ $\forall i\in \left\langle
r\right\rangle ,$ then the average is the fair or uniform allocation of the
resource across the graph; an application (see, \textit{e.g.}, [22, p.
228]): clock synchronization (by averaging) in distributed systems.

\smallskip

For more information concerning the above things, see, \textit{e.g.}, [15],
[27], and [30].

\smallskip

Recall that the DeGroot model and its consensus problem and the consensus
problem in computer science have some things in common with the finite
Markov chains --- some references are given in the first paragraph before
Theorem 2.2 (in Section 2).

\smallskip

\textbf{Remark 3.2.} (a) Concerning the (total) distributed consensuses,
obviously, we have a theorem similar to Theorem 2.2 and a remark similar to
Remark 2.3.

\smallskip

(b) Concerning the finite-time partial and total distributed consensuses
(for reaching a partial or total distributed consensus in a finite time),
only when the matrices at times $t+1,$ $t+2,$ $...$ are (see Definition 1.3)
stochastic real (the matrices at times $1,$ $2,$ $...,$ $t$ being
(stochastic or nonstochastic) real), we have a theorem similar to Theorem
2.5 and remarks similar to Remarks 2.4(c) and 2.6(b)-(c).

\smallskip

(c) Concerning distributed averaging, if we construct/find $t$ matrices, $%
W_{1},$ $W_{2},$ $...,$ $W_{t}$ ($t\geq 1$), such that 
\[
W_{t}W_{t-1}...W_{1}=e^{\prime }\left( \frac{1}{r},\frac{1}{r},...,\frac{1}{r%
}\right) , 
\]
\noindent it makes no sense to construct other matrices --- however, if the
reader wants to construct a sequence of matrices, $\left( W_{n}\right)
_{n\geq 1},$ such that $W_{\infty }=W_{t}W_{t-1}...W_{1},$ she/he can take, 
\textit{e.g.}, $W_{t+1}=W_{t+2}=...=I_{r},$ where $I_{r}$ is the identity
matrix of dimension $r.$

\smallskip

Below we give a result for finite-time (total) distributed averaging --- an
application of the $G$ method (more precisely, of Theorem 1.10)...

\smallskip

\textbf{THEOREM 3.3.}\textit{\ Let }$\mathcal{G}=\left( \mathcal{V},\mathcal{%
E}\right) $\textit{\ be a }(\textit{nondirected simple finite})\textit{\
connected graph.\ Suppose that }$\left| \mathcal{V}\right| =2^{m},$\textit{\ 
}$m\geq 1.$\textit{\ Suppose that }$\mathcal{G}$\textit{\ has a spanning
subgraph isomorphic to the }$m$\textit{-cube graph, }$\mathcal{Q}_{m}$%
\textit{. Then there exist the weight matrices }$W_{1},$ $W_{2},$ $...,$ $%
W_{m}$ \textit{such that distributed averaging is performed for }$\mathcal{G}
$\textit{\ in }$m$ \textit{steps }(\textit{for any initial value vector }$%
q_{0}$). (See Remark 3.2(c) again.)

\smallskip

\textbf{Proof.} It is sufficient to consider the case when $\mathcal{G}=%
\mathcal{Q}_{m}$ --- Theorem 3.3 gives the best possible number of steps, $%
m, $ for distributed averaging for this graph because the diameter of $%
\mathcal{Q}_{m}$ is equal to $m.$ Further, we consider this case, and using
the $G$ method, we show that distributed averaging is performed in $m$ steps.

Consider that $\mathcal{V}=\left\{ 0,1\right\} ^{m}$ and 
\[
\mathcal{E}=\left\{ \left\{ x,y\right\} \left| \text{ }x,y\in \left\{
0,1\right\} ^{m}\text{ and }H\left( x,y\right) =1\right. \right\} , 
\]
\noindent where $H$ is the Hamming distance (the edge set $\mathcal{E}$ is a
set of 2-element subsets of $\left\{ 0,1\right\} ^{m}$; $x=\left(
x_{1},x_{2},...,x_{m}\right) ...$). Consider the stochastic matrices ---
stochastic matrices here! --- $W_{1},W_{2},...,W_{m}$ (see Remark 3.2(c)
again), 
\[
\left( W_{l}\right) _{\left( x_{1},x_{2},...,x_{m}\right) \left(
x_{1},x_{2},...,x_{m}\right) }= 
\]

\[
=\left( W_{l}\right) _{\left( x_{1},x_{2},...,x_{m}\right) \left(
x_{1},x_{2},...,x_{m-l},\left( x_{m-l+1}+1\right) \func{mod}%
2,x_{m-l+2},...,x_{m}\right) }=\frac{1}{2}, 
\]
\noindent $\forall l\in \left\langle m\right\rangle ,$ $\forall \left(
x_{1},x_{2},...,x_{m}\right) \in \left\{ 0,1\right\} ^{m}$ ($\left(
W_{l}\right) _{\left( x_{1},x_{2},...,x_{m}\right) \left(
x_{1},x_{2},...,x_{m}\right) }$ is the entry $\left( \left(
x_{1},x_{2},...,x_{m}\right) ,\left( x_{1},x_{2},...,x_{m}\right) \right) $
of $W_{l}$...; $x_{m-l+2},x_{m-l+3},,...,x_{m}$ vanish when $l=1$; $%
x_{1},x_{2},...,x_{m-l}$ vanish when $l=m$; obviously, $%
W_{1},W_{2},...,W_{m} $ are symmetric matrices). Let 
\[
U_{\left( x_{1},x_{2},...,x_{v}\right) }= 
\]
\[
=\!\left\{ \left( y_{1},y_{2},...,y_{m}\right) \left| \text{ }\left(
y_{1},y_{2},...,y_{m}\right) \!\in \!\left\{ 0,1\right\} ^{m}\text{ and }%
y_{1}=x_{1},y_{2}=x_{2},...,y_{v}=x_{v}\right. \right\} \!\!, 
\]
\noindent $\forall v\in \left\langle m\right\rangle .$ Consider the
partitions 
\[
\Delta _{t}=\left( U_{\left( x_{1},x_{2},...,x_{m-t+1}\right) }\right)
_{\left( x_{1},x_{2},...,x_{m-t+1}\right) \in \left\{ 0,1\right\} ^{m-t+1}},%
\text{ }\forall t\in \left\langle m\right\rangle , 
\]
\noindent and 
\[
\Delta _{m+1}=\left( \left\{ 0,1\right\} ^{m}\right) 
\]
\noindent of $\left\{ 0,1\right\} ^{m}.$ Obviously, 
\[
\Delta _{1}=\left( U_{\left( x_{1},x_{2},...,x_{m}\right) }\right) _{\left(
x_{1},x_{2},...,x_{m}\right) \in \left\{ 0,1\right\} ^{m}}=\left( \left\{
x\right\} \right) _{x\in \left\{ 0,1\right\} ^{m}}. 
\]
\noindent We have $W_{l}\in G_{\Delta _{l+1},\Delta _{l}},\forall l\in
\left\langle m\right\rangle $ --- we prove this statement. Let $l\in
\left\langle m\right\rangle $ and $K\in \Delta _{l+1}.$ It follows that 
\[
K=\left\{ 
\begin{array}{ll}
\left\{ 0,1\right\} ^{m}\smallskip & \text{if }l=m, \\ 
U_{\left( x_{1},x_{2},...,x_{m-l}\right) }\text{ for some }\left(
x_{1},x_{2},...,x_{m-l}\right) \in \left\{ 0,1\right\} ^{m-l} & \text{if }%
l\in \left\langle m-1\right\rangle .%
\end{array}
\right. 
\]

$\smallskip $

\textit{Case} 1. $l=m.$ In this case, $\Delta _{m+1}=\left( \left\{
0,1\right\} ^{m}\right) $ and $\Delta _{m}=\left( U_{\left( 0\right)
},U_{\left( 1\right) }\right) .$ Let $\left( x_{1},x_{2},...,x_{m}\right)
\in \left\{ 0,1\right\} ^{m}.$ Since $\left\{ 0,1\right\} ^{m}=U_{\left(
0\right) }\cup U_{\left( 1\right) },$ it follows that 
\[
\left( x_{1},x_{2},...,x_{m}\right) \in U_{\left( 0\right) }\text{ or }%
\left( x_{1},x_{2},...,x_{m}\right) \in U_{\left( 1\right) }. 
\]
\noindent If $\left( x_{1},x_{2},...,x_{m}\right) \in U_{\left( 0\right) },$
then 
\[
\left( \left( x_{1}+1\right) \func{mod}2,x_{2},...,x_{m}\right) \in
U_{\left( 1\right) }; 
\]
\noindent if $\left( x_{1},x_{2},...,x_{m}\right) \in U_{\left( 1\right) },$
then 
\[
\left( \left( x_{1}+1\right) \func{mod}2,x_{2},...,x_{m}\right) \in
U_{\left( 0\right) }. 
\]
\noindent Further, by the definition of $W_{m},$ we have $W_{m}\in G_{\Delta
_{m+1},\Delta _{m}}.$

$\smallskip $

\textit{Case} 2. $l\in \left\langle m-1\right\rangle .$ In this case, 
\[
K=U_{\left( x_{1},x_{2},...,x_{m-l}\right) }\text{ for some }\left(
x_{1},x_{2},...,x_{m-l}\right) \in \left\{ 0,1\right\} ^{m-l}. 
\]
\noindent Further, we have $K=U_{\left( x_{1},x_{2},...,x_{m-l},0\right)
}\cup U_{\left( x_{1},x_{2},...,x_{m-l},1\right) },$ where, obviously, $%
U_{\left( x_{1},x_{2},...,x_{m-l},0\right) },$ $U_{\left(
x_{1},x_{2},...,x_{m-l},1\right) }\in \Delta _{l}.$ Let $\left(
x_{1},x_{2},...,x_{m}\right) \in K.$ It follows that 
\[
\left( x_{1},x_{2},...,x_{m}\right) \in U_{\left(
x_{1},x_{2},...,x_{m-l},0\right) }\text{ or }\left(
x_{1},x_{2},...,x_{m}\right) \in U_{\left( x_{1},x_{2},...,x_{m-l},1\right)
}. 
\]
\noindent If $\left( x_{1},x_{2},...,x_{m}\right) \in U_{\left(
x_{1},x_{2},...,x_{m-l},0\right) },$ then 
\[
x_{m-l+1}=0\text{ and }\left(
x_{1},x_{2},...,x_{m-l},1,x_{m-l+2},...,x_{m}\right) \in U_{\left(
x_{1},x_{2},...,x_{m-l},1\right) }; 
\]
\noindent if $\left( x_{1},x_{2},...,x_{m}\right) \in U_{\left(
x_{1},x_{2},...,x_{m-l},1\right) },$ then 
\[
x_{m-l+1}=1\text{ and }\left(
x_{1},x_{2},...,x_{m-l},0,x_{m-l+2},...,x_{m}\right) \in U_{\left(
x_{1},x_{2},...,x_{m-l},0\right) }. 
\]
\noindent Further, by the definition of $W_{l},$ we have $W_{l}\in G_{\Delta
_{l+1},\Delta _{l}}.$ The statement was proved. Further, since $W_{l}\in
G_{\Delta _{l+1},\Delta _{l}},$ $\forall l\in \left\langle m\right\rangle ,$
by Theorem 1.10, $W_{m}W_{m-1}...W_{1}$ is a stable matrix and 
\[
W_{m}W_{m-1}...W_{1}=e^{\prime }W_{m}^{-+}W_{m-1}^{-+}...W_{1}^{-+}. 
\]
\noindent Since 
\[
W_{m}^{-+}W_{m-1}^{-+}...W_{1}^{-+}= 
\]
\[
=\left( 
\begin{array}{ll}
\frac{1}{2} & \frac{1}{2}%
\end{array}
\right) \left( 
\begin{array}{cccc}
\frac{1}{2} & \frac{1}{2} &  &  \\ 
&  & \frac{1}{2} & \frac{1}{2}%
\end{array}
\right) ...\left( 
\begin{array}{ccccccc}
\frac{1}{2} & \frac{1}{2} &  &  &  &  &  \\ 
&  & \frac{1}{2} & \frac{1}{2} &  &  &  \\ 
&  &  &  & \ddots &  &  \\ 
&  &  &  &  & \frac{1}{2} & \frac{1}{2}%
\end{array}
\right) = 
\]

\[
=\left( 
\begin{array}{cccc}
\frac{1}{2^{m}} & \frac{1}{2^{m}} & ... & \frac{1}{2^{m}}%
\end{array}
\right) , 
\]
\noindent we have 
\[
q_{m}=\left( W_{m}W_{m-1}...W_{1}q_{0}^{\prime }\right) ^{\prime }=\left( 
\frac{q_{0\left( 0,0,...,0\right) }}{2^{m}},\frac{q_{0\left(
0,0,...,0,1\right) }}{2^{m}},...,\frac{q_{0\left( 1,1,...,1\right) }}{2^{m}}%
\right) 
\]
\noindent ($\mathcal{V}=\left\{ 0,1\right\} ^{m},$ so, $\left| \mathcal{V}%
\right| =2^{m}$; $q_{0\left( 0,0,...,0\right) }$ is the component $\left(
0,0,...,0\right) $ of (initial value vector) $q_{0}$...). So, distributed
averaging is performed in $m$ steps. 
\endproof%

\smallskip

\textbf{Example 3.4.} To illustrate Theorem 3.3, we consider the $3$-cube
graph, $\mathcal{Q}_{3}$. In this case, we have $\mathcal{V}=\left\{
0,1\right\} ^{3}$ and 
\[
\mathcal{E}=\left\{ \left\{ x,y\right\} \left| \text{ }x,y\in \left\{
0,1\right\} ^{3}\text{ and }H\left( x,y\right) =1\right. \right\} = 
\]
\[
=\left\{ \left\{ \left( 0,0,0\right) ,\left( 0,0,1\right) \right\} ,\left\{
\left( 0,0,0\right) ,\left( 0,1,0\right) \right\} ,...,\left\{ \left(
1,1,0\right) ,\left( 1,1,1\right) \right\} \right\} . 
\]
\noindent Further, we have 
\[
\begin{array}{cc}
W_{1}= 
\begin{array}{c}
\left( 0,0,0\right) \smallskip \\ 
\left( 0,0,1\right) \smallskip \\ 
\left( 0,1,0\right) \smallskip \\ 
\left( 0,1,1\right) \smallskip \\ 
\left( 1,0,0\right) \smallskip \\ 
\left( 1,0,1\right) \smallskip \\ 
\left( 1,1,0\right) \smallskip \\ 
\left( 1,1,1\right)%
\end{array}
& \left( 
\begin{array}{cccccccc}
\frac{1}{2}\smallskip & \frac{1}{2} &  &  &  &  &  &  \\ 
\frac{1}{2} & \frac{1}{2}\smallskip &  &  &  &  &  &  \\ 
&  & \frac{1}{2}\smallskip & \frac{1}{2} &  &  &  &  \\ 
&  & \frac{1}{2} & \frac{1}{2}\smallskip &  &  &  &  \\ 
&  &  &  & \frac{1}{2}\smallskip & \frac{1}{2} &  &  \\ 
&  &  &  & \frac{1}{2} & \frac{1}{2}\smallskip &  &  \\ 
&  &  &  &  &  & \frac{1}{2}\smallskip & \frac{1}{2} \\ 
&  &  &  &  &  & \frac{1}{2} & \frac{1}{2}%
\end{array}
\right)%
\end{array}
\]
\noindent \noindent (the columns are labeled similarly, \textit{i.e.}, using
the lexicographic order from left to right), 
\[
\begin{array}{cc}
W_{2}= 
\begin{array}{c}
\left( 0,0,0\right) \smallskip \\ 
\left( 0,0,1\right) \smallskip \\ 
\left( 0,1,0\right) \smallskip \\ 
\left( 0,1,1\right) \smallskip \\ 
\left( 1,0,0\right) \smallskip \\ 
\left( 1,0,1\right) \smallskip \\ 
\left( 1,1,0\right) \smallskip \\ 
\left( 1,1,1\right)%
\end{array}
& \left( 
\begin{array}{cccccccc}
\frac{1}{2}\smallskip & 0 & \frac{1}{2} & 0 &  &  &  &  \\ 
0 & \frac{1}{2}\smallskip & 0 & \frac{1}{2} &  &  &  &  \\ 
\frac{1}{2} & 0 & \frac{1}{2}\smallskip & 0 &  &  &  &  \\ 
0 & \frac{1}{2} & 0 & \frac{1}{2}\smallskip &  &  &  &  \\ 
&  &  &  & \frac{1}{2}\smallskip & 0 & \frac{1}{2} & 0 \\ 
&  &  &  & 0 & \frac{1}{2}\smallskip & 0 & \frac{1}{2} \\ 
&  &  &  & \frac{1}{2} & 0 & \frac{1}{2}\smallskip & 0 \\ 
&  &  &  & 0 & \frac{1}{2} & 0 & \frac{1}{2}%
\end{array}
\right) ,%
\end{array}
\]
\[
\begin{array}{cc}
W_{3}= 
\begin{array}{c}
\left( 0,0,0\right) \smallskip \\ 
\left( 0,0,1\right) \smallskip \\ 
\left( 0,1,0\right) \smallskip \\ 
\left( 0,1,1\right) \smallskip \\ 
\left( 1,0,0\right) \smallskip \\ 
\left( 1,0,1\right) \smallskip \\ 
\left( 1,1,0\right) \smallskip \\ 
\left( 1,1,1\right)%
\end{array}
& \left( 
\begin{array}{cccccccc}
\frac{1}{2}\smallskip & 0 & 0 & 0 & \frac{1}{2} & 0 & 0 & 0 \\ 
0 & \frac{1}{2}\smallskip & 0 & 0 & 0 & \frac{1}{2} & 0 & 0 \\ 
0 & 0 & \frac{1}{2}\smallskip & 0 & 0 & 0 & \frac{1}{2} & 0 \\ 
0 & 0 & 0 & \frac{1}{2}\smallskip & 0 & 0 & 0 & \frac{1}{2} \\ 
\frac{1}{2} & 0 & 0 & 0 & \frac{1}{2}\smallskip & 0 & 0 & 0 \\ 
0 & \frac{1}{2} & 0 & 0 & 0 & \frac{1}{2}\smallskip & 0 & 0 \\ 
0 & 0 & \frac{1}{2} & 0 & 0 & 0 & \frac{1}{2}\smallskip & 0 \\ 
0 & 0 & 0 & \frac{1}{2} & 0 & 0 & 0 & \frac{1}{2}%
\end{array}
\right) .%
\end{array}
\]
\noindent Further, we have 
\[
U_{\left( x_{1},x_{2},x_{3}\right) }=\left\{ \left( x_{1},x_{2},x_{3}\right)
\right\} ,\text{ }\forall \left( x_{1},x_{2},x_{3}\right) \in \left\{
0,1\right\} ^{3}, 
\]
\[
U_{\left( 0,0\right) }=\left\{ \left( 0,0,0\right) ,\left( 0,0,1\right)
\right\} ,U_{\left( 0,1\right) }=\left\{ \left( 0,1,0\right) ,\left(
0,1,1\right) \right\} , 
\]
\[
U_{\left( 1,0\right) }=\left\{ \left( 1,0,0\right) ,\left( 1,0,1\right)
\right\} ,U_{\left( 1,1\right) }=\left\{ \left( 1,1,0\right) ,\left(
1,1,1\right) \right\} , 
\]
\[
U_{\left( 0\right) }=U_{\left( 0,0\right) }\cup U_{\left( 0,1\right)
}=\left\{ \left( 0,0,0\right) ,\left( 0,0,1\right) ,\left( 0,1,0\right)
,\left( 0,1,1\right) \right\} , 
\]
\[
U_{\left( 1\right) }=U_{\left( 1,0\right) }\cup U_{\left( 1,1\right)
}=\left\{ \left( 1,0,0\right) ,\left( 1,0,1\right) ,\left( 1,1,0\right)
,\left( 1,1,1\right) \right\} . 
\]
\noindent Further, we have 
\[
\Delta _{1}=\left( \left\{ \left( x_{1},x_{2},x_{3}\right) \right\} \right)
_{\left( x_{1},x_{2},x_{3}\right) \in \left\{ 0,1\right\} ^{3}}=\left(
\left\{ \left( 0,0,0\right) \right\} ,\left\{ \left( 0,0,1\right) \right\}
,...,\left\{ \left( 1,1,1\right) \right\} \right) , 
\]
\[
\Delta _{2}=\left( U_{\left( 0,0\right) },U_{\left( 0,1\right) },U_{\left(
1,0\right) },U_{\left( 1,1\right) }\right) , 
\]
\[
\Delta _{3}=\left( U_{\left( 0\right) },U_{\left( 1\right) }\right) , 
\]
\[
\Delta _{4}=\left( \left\{ 0,1\right\} ^{3}\right) 
\]
\noindent ($\Delta _{4}$ is the improper (degenerate) partition of $\left\{
0,1\right\} ^{3}$ --- $\Delta _{4}$ has only one set). Further, we have $%
W_{1}\in G_{\Delta _{2},\Delta _{1}},$ $W_{2}\in G_{\Delta _{3},\Delta
_{2}}, $ $W_{3}\in G_{\Delta _{4},\Delta _{3}};$ further, we have 
\[
W_{3}^{-+}W_{2}^{-+}W_{1}^{-+}= 
\]
\[
=\left( 
\begin{array}{cc}
\frac{1}{2} & \frac{1}{2}%
\end{array}
\right) \left( 
\begin{array}{cccc}
\frac{1}{2} & \frac{1}{2} &  &  \\ 
&  & \frac{1}{2} & \frac{1}{2}%
\end{array}
\right) \left( 
\begin{array}{cccccccc}
\frac{1}{2} & \frac{1}{2} &  &  &  &  &  &  \\ 
&  & \frac{1}{2} & \frac{1}{2} &  &  &  &  \\ 
&  &  &  & \frac{1}{2} & \frac{1}{2} &  &  \\ 
&  &  &  &  &  & \frac{1}{2} & \frac{1}{2}%
\end{array}
\right) = 
\]
\[
=\left( 
\begin{array}{cccc}
\frac{1}{2^{3}} & \frac{1}{2^{3}} & \cdots & \frac{1}{2^{3}}%
\end{array}
\right) . 
\]
\noindent By Theorem 1.10 ($W_{3}W_{2}W_{1}=e^{\prime
}W_{3}^{-+}W_{2}^{-+}W_{1}^{-+}$) or by direct computation, 
\[
W_{3}W_{2}W_{1}=\left( 
\begin{array}{cccc}
\frac{1}{2^{3}}\smallskip & \frac{1}{2^{3}} & \cdots & \frac{1}{2^{3}} \\ 
\frac{1}{2^{3}} & \frac{1}{2^{3}}\smallskip & \cdots & \frac{1}{2^{3}} \\ 
\vdots & \vdots & \cdots \smallskip & \vdots \\ 
\frac{1}{2^{3}} & \frac{1}{2^{3}} & \cdots & \frac{1}{2^{3}}%
\end{array}
\right) . 
\]
\noindent Finally, we have ($\mathcal{V}=\left\{ 0,1\right\} ^{3},$ so, $%
\left| \mathcal{V}\right| =2^{3}$; $q_{0\left( 0,0,0\right) }$ is the
component $\left( 0,0,0\right) $ of $q_{0}$...) 
\[
q_{3}=\left( W_{3}W_{2}W_{1}q_{0}^{\prime }\right) ^{\prime }=\left( \frac{%
q_{0\left( 0,0,0\right) }}{2^{3}},\frac{q_{0\left( 0,0,1\right) }}{2^{3}}%
,...,\frac{q_{0\left( 1,1,1\right) }}{2^{3}}\right) 
\]
\noindent --- distributed averaging is performed at time/step $3$ (for any
initial value vector $q_{0}$). To understand Theorem 3.3 and this example
better, the reader, if she/he wants, can compute the value vectors $q_{1}$
and $q_{2}$. Note that the matrices $W_{1},W_{2},$ and $W_{3}$ are
symmetric, but by Theorem 1.16, $W_{2}$ and $W_{3}$ can be replaced with
nonsymmetric ones.

\smallskip

\textbf{Remark 3.5.} (Research work.) Theorem 3.3 can be extended for any
connected graph with $n_{1}n_{2}...n_{t}$ vertices, where $n_{1},$ $n_{2},$ $%
...,$ $n_{t}\geq 2$ (preferably prime numbers) and having a spanning
subgraph for which Theorem 1.10 can be applied for $t$ stochastic
matrices... --- we will perform distributed averaging in $t$ steps.

\smallskip

For the above remark, we consider the next example.

\smallskip

\textbf{Example 3.6.} (Suggested by Example 2.11 from [17], which refers to
the uniform generation of permutations of order $n$.) Consider the graph $%
\mathcal{G}_{3}=\left( \mathcal{V}_{3},\mathcal{E}_{3}\right) ,$ where $%
\mathcal{V}_{3}=\mathbb{S}_{3},$ $\mathbb{S}_{3}=$ the set of permutations
of order $3,$ and $\mathcal{E}_{3}=$ the union of edge sets of nondirected
graphs of matrices $W_{1}$ and $W_{2}$, where 
\[
\begin{array}{cc}
W_{1}= 
\begin{array}{c}
\left( 123\right) \smallskip \\ 
\left( 132\right) \smallskip \\ 
\left( 213\right) \smallskip \\ 
\left( 231\right) \smallskip \\ 
\left( 321\right) \smallskip \\ 
\left( 312\right)%
\end{array}
& \left( 
\begin{array}{cccccc}
\frac{1}{2}\smallskip & \frac{1}{2} &  &  &  &  \\ 
\frac{1}{2} & \frac{1}{2}\smallskip &  &  &  &  \\ 
&  & \frac{1}{2}\smallskip & \frac{1}{2} &  &  \\ 
&  & \frac{1}{2} & \frac{1}{2}\smallskip &  &  \\ 
&  &  &  & \frac{1}{2}\smallskip & \frac{1}{2} \\ 
&  &  &  & \frac{1}{2} & \frac{1}{2}%
\end{array}
\right) ,%
\end{array}
\]
\[
\begin{array}{cc}
W_{2}= 
\begin{array}{c}
\left( 123\right) \smallskip \\ 
\left( 132\right) \smallskip \\ 
\left( 213\right) \smallskip \\ 
\left( 231\right) \smallskip \\ 
\left( 321\right) \smallskip \\ 
\left( 312\right)%
\end{array}
& \left( 
\begin{array}{cccccc}
\frac{1}{3}\smallskip & 0 & \frac{1}{3} & 0 & \frac{1}{3} & 0 \\ 
0 & \frac{1}{3}\smallskip & 0 & \frac{1}{3} & 0 & \frac{1}{3} \\ 
\frac{1}{3} & 0 & \frac{1}{3}\smallskip & 0 & 0 & \frac{1}{3} \\ 
0 & \frac{1}{3} & 0 & \frac{1}{3}\smallskip & \frac{1}{3} & 0 \\ 
\frac{1}{3} & 0 & 0 & \frac{1}{3} & \frac{1}{3}\smallskip & 0 \\ 
0 & \frac{1}{3} & \frac{1}{3} & 0 & 0 & \frac{1}{3}%
\end{array}
\right) .%
\end{array}
\]
\noindent The matrices $W_{1}$ and $W_{2}$ were obtained by the swapping
method --- for this method, see, \textit{e.g.}, [8, pp. 645$-$646]. By
Theorem 1.10 ($W_{1}\in G_{\Delta _{2},\Delta _{1}},$ $W_{2}\in G_{\Delta
_{3},\Delta _{2}},$ where...) or by direct computation, 
\[
W_{2}W_{1}=\left( 
\begin{array}{cccc}
\frac{1}{6}\smallskip & \frac{1}{6} & \cdots & \frac{1}{6} \\ 
\frac{1}{6} & \frac{1}{6}\smallskip & \cdots & \frac{1}{6} \\ 
\vdots & \vdots & \cdots \smallskip & \vdots \\ 
\frac{1}{6} & \frac{1}{6} & \cdots & \frac{1}{6}%
\end{array}
\right) , 
\]
\noindent so, distributed averaging is performed, in this case, at time $2$ (%
$3!=6=2\cdot 3$). Note that the graph $\mathcal{G}_{3}$ is a $3$-regular
graph and is not isomorphic to the triangular prism graph, which is also a $%
3 $-regular graph --- $\mathcal{G}_{3}$ can be replaced with any graph with $%
3! $ vertices which has a spanning subgraph isomorphic to $\mathcal{G}_{3}.$
For the triangular prism graph, using the $G$ method, distributed averaging
can also be performed at time 2 --- an exercise for the reader.

\smallskip

The above example can be generalized for $\mathcal{G}_{n}=\left( \mathcal{V}%
_{n},\mathcal{E}_{n}\right) ,$ $\mathcal{V}_{n}=\mathbb{S}_{n},$ $\mathbb{S}%
_{n}=$ the set of permutations of order $n$,... Due to the swapping method,
this graph must be an $\frac{n\left( n-1\right) }{2}$-regular graph --- a
necessary but not sufficient condition (see the above example).

\smallskip

\textbf{Remark 3.7.} (a) For the DeGroot model on distributed systems,
labelling the vertices of (involved) graphs is an important problem --- a
good/

\noindent suitable labelling can lead to a good/fast algorithm (see Theorem
3.3 and its proof and Examples 3.4 and 3.6).

\smallskip

(b) Theorem 3.3 and its proof and Example 3.6 (see also Remark 3.5) lead to
possible cases for the DeGroot submodel (from Section 2) on graphs/networks,
in each case, a consensus being reached in a finite time --- for the DeGroot
model on graphs, see, \textit{e.g.}, [9] and [12].

\smallskip

For the distributed consensuses (in particular, for distributed averaging),
we could combine the $G$ method with other methods, obtaining certain hybrid
methods: we could use the $G$ method for subgraphs, we could use the $G$
method together with the flooding (for this method, see, \textit{e.g.} [30,
p. 65], see, \textit{e.g.}, also [22, p. 228]), we could use a leader vertex
or more... This idea can be developed --- another research work ---, but we
give only an example, the next example.

\smallskip

\textbf{Example 3.8.} Consider the graph from Fig. 1 in [30] (this graph is
also considered in [27, Fig. 1]) --- the weights from there are not taken
into account. This graph, say, $\mathcal{G},$ can be obtained from the wheel
graph $\mathcal{W}_{5}$ and complete graph $\mathcal{K}_{5}$ by edge
merging. Consider that the vertices of $\mathcal{K}_{5}$ are $1,$ $2,$ $3,$ $%
4,$ $5$ and the vertices of $\mathcal{W}_{5}$ are $1,$ $2,$ $6,$ $7,$ $8$,
and that $8$ is the universal vertex of $\mathcal{W}_{5},$ $6$ is adjacent
to $1$ ($6$ is also adjacent to $8$), and $7$ is adjacent to $2$ ($7$ is
also adjacent to $6$ and $8$)--- $\left\{ 1,2\right\} $ is the common edge
of $\mathcal{W}_{5}$ and $\mathcal{K}_{5}.$ Consider that the initial value
vector is $q_{0}.$ First, we consider the subgraph $\mathcal{W}_{5}$ of $%
\mathcal{G}.$ For this subgraph, we consider that the initial value vector
is $q_{0}^{\left( 1\right) }=q_{0}\mid _{\left\{ 1,2,6,7,8\right\} },$ $%
q_{0}\mid _{\left\{ 1,2,6,7,8\right\} }=$ the restriction of $q_{0}$ to $%
\left\{ 1,2,6,7,8\right\} .$ Consider the matrix $W_{1}^{\left( 1\right) }$
with rows and columns $1,$ $2,$ $6,$ $7,$ $8,$%
\[
W_{1}^{\left( 1\right) }=\left( 
\begin{array}{cccc}
\frac{1}{5}\smallskip & \frac{1}{5} & \cdots & \frac{1}{5} \\ 
\frac{1}{5} & \frac{1}{5}\smallskip & \cdots & \frac{1}{5} \\ 
\vdots & \vdots & \cdots \smallskip & \vdots \\ 
\frac{1}{5} & \frac{1}{5} & \cdots & \frac{1}{5}%
\end{array}
\right) . 
\]
\noindent $W_{1}^{\left( 1\right) }$ is a stable matrix, and, therefore, $%
W_{1}^{\left( 1\right) }\in G_{\Delta _{2}^{\left( 1\right) },\Delta
_{1}^{\left( 1\right) }},$ where $\Delta _{1}^{\left( 1\right) }=\left(
\left\{ 1\right\} ,\left\{ 2\right\} ,\left\{ 6\right\} ,\left\{ 7\right\}
,\left\{ 8\right\} \right) $ and $\Delta _{2}^{\left( 1\right) }=\left(
\left\{ 1,2,6,7,8\right\} \right) $ --- a trivial case for the $G$ method,
more exactly, for Theorem 1.10. Further, from 
\[
\left( q_{1}^{\left( 1\right) }\right) ^{\prime }=W_{1}^{\left( 1\right)
}\left( q_{0}^{\left( 1\right) }\right) ^{\prime }=W_{1}^{\left( 1\right)
}\left( q_{0}\mid _{\left\{ 1,2,6,7,8\right\} }\right) ^{\prime }, 
\]
\noindent we obtain 
\[
q_{1}^{\left( 1\right) }=\left( \frac{1}{5}\sum\limits_{i\in \left\{
1,2,6,7,8\right\} }q_{0i},\frac{1}{5}\sum\limits_{i\in \left\{
1,2,6,7,8\right\} }q_{0i},...,\frac{1}{5}\sum\limits_{i\in \left\{
1,2,6,7,8\right\} }q_{0i}\right) , 
\]
\noindent $q_{1}^{\left( 1\right) }$ is the value vector at time 1 for the
subgraph $\mathcal{W}_{5}.$ Second, we consider the subgraph of $\mathcal{G}$
with vertices $1,$ $3,$ $4,$ $5$ and isomorphic to the complete graph $%
\mathcal{K}_{4}.$ For this subgraph, we consider that the initial value
vector is 
\[
\left( \frac{1}{5}\sum\limits_{i\in \left\{ 1,2,6,7,8\right\}
}q_{0i},q_{03},q_{04},q_{05}\right) :=q_{0}^{\left( 2\right) }. 
\]
\noindent Consider the matrix $W_{1}^{\left( 2\right) }$ with rows and
columns $1,$ $3,$ $4,$ $5,$%
\[
W_{1}^{\left( 2\right) }=\left( 
\begin{array}{cccc}
\frac{5}{8}\smallskip & \frac{1}{8} & \frac{1}{8} & \frac{1}{8} \\ 
\frac{5}{8} & \frac{1}{8}\smallskip & \frac{1}{8} & \frac{1}{8} \\ 
\frac{5}{8} & \frac{1}{8} & \frac{1}{8}\smallskip & \frac{1}{8} \\ 
\frac{5}{8} & \frac{1}{8} & \frac{1}{8} & \frac{1}{8}%
\end{array}
\right) . 
\]
\noindent $W_{1}^{\left( 2\right) }$ is also a stable matrix. Further, from 
\[
\left( q_{1}^{\left( 2\right) }\right) ^{\prime }=W_{1}^{\left( 2\right)
}\left( q_{0}^{\left( 2\right) }\right) ^{\prime }, 
\]
\noindent we obtain 
\[
q_{1}^{\left( 2\right) }=\left( \frac{1}{8}\sum\limits_{i\in \left\langle
8\right\rangle }q_{0i},\text{ }\frac{1}{8}\sum\limits_{i\in \left\langle
8\right\rangle }q_{0i},\text{ }\frac{1}{8}\sum\limits_{i\in \left\langle
8\right\rangle }q_{0i},\text{ }\frac{1}{8}\sum\limits_{i\in \left\langle
8\right\rangle }q_{0i}\right) , 
\]
\noindent $q_{1}^{\left( 2\right) }$ is the value vector at time 1 for the
above subgraph. Consider that the vertex $1$ (it can be a leader) transmits
the average $\frac{1}{8}\sum\limits_{i\in \left\langle 8\right\rangle
}q_{0i} $ to the vertices $2,$ $6,$ and $8$ and, further, that the vertex $2$
(this can also be a leader) or $6$ transmits the above average to $7$ ---
and thus distributed averaging is performed. Note that we can use $\mathcal{K%
}_{5}$ instead of the subgraph of $\mathcal{G}$ with vertices $1,$ $3,$ $4,$ 
$5$ and isomorphic to $\mathcal{K}_{4}...$; finally, the vertices $1$ and $2$
(these can be leaders) will transmit the average $\frac{1}{8}%
\sum\limits_{i\in \left\langle 8\right\rangle }q_{0i}$ to the vertices $6,$ $%
7,$ and $8.$ The study of this alternative way is left to the reader.

\smallskip

For the (nondirected simple finite) connected graphs with at least two
vertices, distributed averaging based on the DeGroot model on distributed
systems and the uniform generation based on Markov chains of the vertices of
a graph taking its edge set into account (when $\left\{ i,j\right\} $ ($%
i\neq j$) is not a edge of the graph, any transition matrix we use has the
entries $\left( i,j\right) $ and $\left( j,i\right) $ equal to 0) are
related problems/topics (... the stochastic matrices from Example 3.4 can be
used to construct a Markov chain for the uniform generation of vertices of
3-cube graph...) 
\[
\text{--- \textit{which of them is harder}?\textit{...}} 
\]

\smallskip

In this article, the framework for the DeGroot model and that for the
DeGroot model on distributed systems we used are general, \textit{i.e.},
homogeneous (when we use just one matrix at any time/step) and
nonhomogeneous (when we use at least two different matrices at least two
different times) 
\[
\text{ --- \textit{which of them are better}?} 
\]
\noindent We believe that the nonhomogeneous frameworks are better --- these
frameworks are supported, among other things, by the DeGroot submodel,
Theorem 3.3 (Example 3.4 is for $m=3$), and Example 3.6 (we used
nonhomogeneous products of matrices in these three places). For the
homogeneous frameworks, see, \textit{e.g.}, [7], [27], and [30].

\bigskip

\begin{center}
\textbf{REFERENCES}
\end{center}

\bigskip \ 

[1] H. Attiya and J. Welch, \textit{Distributed Computing}: \textit{%
Fundamentals, Simulations and Advances Topics}, 2nd Edition. Wiley, Hoboken,
NJ, 2004.

[2] R.L. Berger, \textit{A necessary and sufficient condition for reaching a
consensus using DeGroot's method}. J. Amer. Statist. Assoc. \textbf{76}
(1981), 415$-$418.

[3] D.P. Bertsekas and J.N. Tsitsiklis, \textit{Parallel and Distributed
Computation}:\textit{\ Numerical Methods}. Prentice-Hall, Englewood Cliff,
NJ, 1989.

[4] S.\ Boccaletti, P. De Lellis, C.I. del Genio, K. Alfaro-Bittner, R.
Criado, S. Jalan, and M. Romance, \textit{The structure and dynamics of
networks with higher order interactions}. Phys. Rep. \textbf{1018} (2023), 1$%
-$64.

[5] S. Chatterjee and E. Seneta, \textit{Towards consensus}:\textit{\ some
convergence theorems on repeated averaging}. J. Appl. Probab. \textbf{14}
(1977), 89$-$97.

[6] G. Coulouris, J. Dollimore, T. Kindberg, and G. Blair, \textit{%
Distributed Systems}:\textit{\ Concepts and Design}, 5th Edition.
Addison-Wesley, Boston, 2012.

[7] M.H. DeGroot, \textit{Reaching a consensus}. J. Amer. Statist. Assoc. 
\textbf{69} (1974), 118$-$121.

[8] L. Devroye, \textit{Non-Uniform Random Variate Generation}.
Springer-Verlag, New York, 1986.

[9] B. Golub and M.O. Jackson, \textit{Na\"{\i}ve learning in social
networks and the wisdom of crowds}. AEJ: Microeconomics \textbf{2} (2010),
112$-$149.

[10] M. Iosifescu, \textit{Finite Markov Processes and Their Applications}.
Wiley, Chichester \& Ed. Tehnic\u{a}, Bucharest, 1980; corrected
republication by Dover, Mineola, N.Y., 2007.

[11] D.L. Isaacson and R.W. Madsen, \textit{Markov Chains}:\textit{\ Theory
and Applications.} Wiley, New York, 1976; republication by Krieger, 1985.

[12] M.O. Jackson, \textit{Social and Economic Networks}. Princeton
University Press, Princeton, 2008.

[13] D.A. Levin and Y. Peres, \textit{Markov Chains and Mixing Times}, 2nd
Edition. AMS, Providence, RI, 2017. With contributions by E.L. Wilmer. With
a chapter on ``Coupling from the Past'' by J.G. Propp and D.B. Wilson.

[14] S. Li, H. Du, and X. Lin, \textit{Finite-time consensus algorithm for
multi-agent systems with double-integrator dynamics}. Automatica J. IFAC 
\textbf{47} (2011), 1706$-$1712.

[15] N.A. Lynch, \textit{Distributed Algorithms}. Morgan Kaufmann
Publishers, Inc., 1996.

[16] U. P\u{a}un, \textit{General }$\Delta $\textit{-ergodic theory, with
some results on simulated annealing}. Math. Rep. (Bucur.) \textbf{13}(%
\textbf{63}) (2011), 171$-$196.

[17] U. P\u{a}un, $G_{\Delta _{1},\Delta _{2}}$ \textit{in action}. Rev.
Roumaine Math. Pures Appl.\textbf{\ 55 }(2010), 387$-$406.

[18] U. P\u{a}un, \textit{A hybrid Metropolis-Hastings chain}. Rev. Roumaine
Math. Pures Appl.\textbf{\ 56 }(2011), 207$-$228.

[19] U. P\u{a}un, $G$\textit{\ method in action}: \textit{from exact
sampling to approximate one}. Rev. Roumaine Math. Pures Appl.\textbf{\ 62 }%
(2017), 413$-$452.

[20] U. P\u{a}un, \textit{Ewens distribution on }$\mathbb{S}_{n}$ \textit{is
a wavy probability distribution with respect to }$n$\textit{\ partitions.}
An. Univ. Craiova Ser. Mat. Inform. \textbf{47} (2020), 1$-$24.

[21] U. P\u{a}un, $G$\textit{\ method in action}: \textit{pivot}$^{\text{+}} 
$\textit{\ algorithm for self-avoiding walk}. ROMAI\ J. \textbf{20} (2024),
no. 2, 81$-$98. With preprint: arXiv:2310.07564.

[22] M. Pease, R. Shostak, and L. Lamport, \textit{Reaching agreement in the
presence of faults}. J. ACM \textbf{27} (1980), 228$-$234.

[23] W. Ren, R.W. Beard, and E.M. Atkins, \textit{A survey of consensus
problems in multi-agent coordination}. In: Proceedings of the American
Control Conference, 2005, pp. 1859$-$1864.

[24] L. Saloff-Coste, \textit{Lectures on finite Markov chains}. In: E. Gin%
\'{e}, G.R. Grimmett, and L. Saloff-Coste (Eds.), Lectures on Probability
Theory and Statistics, Lecture Notes in Mathematics 1665, Springer, Berlin,
1997.

[25] L. Saloff-Coste and J. Z\'{u}\~{n}iga, \textit{Merging and stability
for time inhomogeneous finite Markov chains}. In: J. Blath, P. Imkeller, and
S. R\oe lly (Eds.), Surveys in Stochastic Processes, 127$-$151, EMS Ser.
Congr. Rep. Eur. Math. Soc., Z\"{u}rich, 2011.

[26] E. Seneta, \textit{Non-negative Matrices and Markov Chains}, 2nd
Edition. Springer-Verlag, Berlin, 1981; revised printing, 2006.

[27] S. Sundaram and C.N. Hadjicostis, \textit{Finite-time distributed
consensus in graphs with time-invariant topologies}. 2007 American Control
Conference, New York, NY, USA, 2007, pp. 711$-$716.

[28] D.D. \v{S}iljak, \textit{Large-Scale Dynamic Systems}: \textit{%
Stability and Structure}. North-Holland, New York, 1978.

[29] W. Xia, J. Liu, M. Cao, K.H. Johansson, and T. Ba\c{s}ar, \textit{%
Generalized Sarymsakov matrices}. IEEE Trans. Automat. Control \textbf{64}
(2019), 3085$-$3100.

[30] L. Xiao and S. Boyd, \textit{Fast linear iterations for distributed
averaging}. Systems Control Lett. \textbf{53} (2004) 65$-$78.

\bigskip\ 
\[
\begin{array}{ccc}
\mathit{October}\text{ }\mathit{20,}\text{ }\mathit{2025} &  & \mathit{%
Romanian\ Academy} \\ 
&  & \mathit{Gheorghe\ Mihoc-Caius\ Iacob}\text{ }\mathit{Institute} \\ 
&  & \mathit{of\ Mathematical\ Statistics}\text{ }\mathit{and\ Applied\
Mathematics} \\ 
&  & \mathit{Calea\ 13\ Septembrie\ nr.\ 13} \\ 
&  & \mathit{050711\ Bucharest\ 5,\ Romania} \\ 
&  & \mathit{upterra@gmail.com}%
\end{array}%
\]

\end{document}